\magnification=1100
\input amstex
\documentstyle{amsppt}
\NoBlackBoxes

\hsize=14cm
\vsize=19cm

\def\a{\alpha}

\topmatter
\title Counting Congruence Subgroups \endtitle
\author Dorian Goldfeld \\ Alexander Lubotzky \\ L\'aszl\'o Pyber
\endauthor
\address D. Goldfeld, Department of Mathematics, Columbia University, NY, NY 10027, USA
\endaddress
\email goldfeld\@columbia.edu \endemail
\address A. Lubotzky, Institute of Mathematics, Hebrew University, Jerusalem 91904, Israel
\endaddress
\email alexlub\@math.huji.ac.il \endemail
\address L. Pyber, A. Renyi Institute of Mathematics,
Hungarian Academy of Sciences, \break
\phantom{x}\hskip 59pt P.O. Box 127,
H-1364 Budapest
,Hungary
\endaddress
\email pyber\@renyi.hu \endemail
\thanks  The first two authors research is  supported in part by the NSF. The third author's Research is supported  in
part by OTKA T 034878. All three authors would like to thank Yale University for its hospitality.
\endthanks

\abstract   Let $\Gamma$ denote the modular group $SL(2,\Bbb Z)$
and $C_n(\Gamma)$ the number of congruence subgroups of $\Gamma$
of index at most $n$. We prove that  $\lim\limits_{n\to \infty}
\frac{\log C_n(\Gamma)}{ (\log n)^2/\log\log n} =
\frac{3-2\sqrt{2}}{4}.$ We also present a very general conjecture
giving an asymptotic estimate for $C_n(\Gamma)$ for general
arithmetic groups. The lower bound of the conjecture is proved
modulo the generalized Riemann hypothesis for Artin-Hecke L-functions,
and in many cases is also proved unconditionally. The upper bound
is proved in full in \cite{LN}.
\endabstract

\endtopmatter

\document
 \baselineskip 15pt

{\bf \S0.  Introduction}
\vskip .10in

Let $k$ be an algebraic number field, $\Cal O$ its ring of
integers, $S$ a
 finite set of valuations of $k$ (containing  all the archimedean
ones), and $\Cal O_S = \big\{x\in k \; \big\vert \; v(x)\geq 0, \;
\forall v\not\in S\big\}$.  Let $G$ be a semisimple, simply
connected, connected algebraic group defined over $k$ with a fixed
embedding into $GL_d$.  Let $\Gamma = G(\Cal O_S) = G\cap
GL_d(\Cal O_S)$ be the corresponding $S$-arithmetic group.  We
assume that $\Gamma$ is an infinite group (equivalently,
$\prod_{\nu\in S} G(k_\nu)$ is not compact).

For every non-zero ideal $I$ of $\Cal O_S$ let $$\Gamma (I) =
\text{Ker}\big(\Gamma\to GL_d(\Cal O_S/I)\big).$$  A subgroup of
$\Gamma$ is called a congruence subgroup if it contains $\Gamma
(I)$ for some $I$.

The topic of counting congruence subgroups has a long history.
Classically, congruence subgroups of the modular group were counted as a function of
the genus of the associated Riemann surface. It was conjectured by
Rademacher that there are only finitely many congruence subgroups
of $SL_2(\Bbb Z)$ of genus zero.  Petersson \cite{Pe, 1974} proved
that the number of all subgroups of index $n$ and fixed genus goes
to infinity exponentially as $n \to \infty.$ Dennin \cite{De,
1975} proved that there are only finitely many congruence
subgroups of  $SL_2(\Bbb Z)$ of given fixed genus and solved
Rademacher's conjecture. It does not seem possible, however, to
accurately count all congruence subgroups of index at most $n$ in
$SL_2(\Bbb Z)$ by using the theory of Riemann surfaces of fixed
genus.

Following  \cite{Lu}, we count congruence subgroups as a function
of the index. For $n > 0$, define
$$C_n(\Gamma ) = \# \big\{\text{congruence subgroups
of $\Gamma$ of index at most $n$}\big\}.$$

\proclaim{Theorem 1}  There exist two positive real numbers
$\alpha_-(\Gamma)$ and $\alpha_+(\Gamma)$ such that for all sufficiently large
positive integers $n$
$$n^{\frac{\log n}{\log\log n}\alpha_{-}} \; \leq  \; C_n(\Gamma) \;
\leq \; n^{\frac{\log n}{\log\log n}\alpha_{+}} .$$
\endproclaim

This theorem is proved in \cite{Lu}, although the proof of the
lower bound presented there requires the prime number theorem on
arithmetic progressions in an interval where its validity depends
on the GRH (generalized Riemann hypothesis for Dirichlet L-functions). By a slight modification of the proof and by
appealing to a theorem of Linnik \cite{Li1, Li2} on the least
prime in an arithmetic progression, the proof can be made
unconditional. Such an approach gives, however, poor estimates for
the constants.

 Following \cite{Lu} we define:
$$\alpha_+(\Gamma) =
\overline{\lim} \; \frac{\log C_{n}(\Gamma )}{\lambda (n)}, \qquad
\alpha_-(\Gamma) = \underline{\lim}  \; \frac{\log
C_{n}(\Gamma)}{\lambda (n)},
$$
where $\lambda (n) = \frac{(\log n)^{2}}{\log\log n}$.
\vskip 8pt
It is not difficult to see that $\alpha_+$ and $\alpha_-$ are
independent of both the choice of the representation of $G$ as a
matrix group and of the choice of $S$.  Hence
$\alpha_\pm$ depend only on $G$ and $k$.  The question whether
$\alpha_+(\Gamma) = \alpha_-(\Gamma)$  and the challenge to
evaluate them   for $\Gamma = SL_2(\Bbb Z)$ and other groups was
presented in \cite{Lu}. Here we prove:

\proclaim{Theorem 2}  We have $\alpha_+(SL_2(\Bbb Z)) =
\alpha_- (SL_2(\Bbb Z)) = \frac{3 - 2\sqrt{2}}{4} =
0.0428932\ldots$
\endproclaim

The proof of the lower bound in Theorem 2 is based on the
Bombieri-Vinogradov Theorem \cite{Bo}, \cite{Da}, \cite{Vi},
i.e., {\it the Riemann hypothesis on the average.}  The upper
bound, on the other hand, is proved by first reducing the problem to a
counting problem for subgroups of abelian groups and then
solving  that extremal counting problem.

In the case of a number field, we will, in fact, show a more remarkable result: the answer is
independent of $\Cal O$! Here, we require the GRH (generalized Riemann hypothesis)
\cite{W} for Hecke and Artin L-functions which states that all non-trivial zeros of such L-functions lie on the critical line.

\proclaim{Theorem 3}  Let $k$ be a number field  with ring of
integers $\Cal O$. Let $S$ be a finite set of primes, and $\Cal
O_S$ as above. Assume GRH  for $k$
and all cyclotomic extensions  $k(\zeta_\ell)$ with $\ell$ a
rational prime and $\zeta_\ell$ a primitive $\ell^{\text{\rm th}}$
root of unity. Then
$$
\alpha_+ (SL_2(\Cal O_S)) = \alpha_- (SL_2(\Cal O_S)) =
\frac{3 - 2\sqrt{2}}{4}.
$$
\endproclaim

The GRH is needed only for establishing the lower bound. It can be
dropped in many cases by appealing to a theorem of Murty and Murty
\cite{MM} which generalizes the Bombieri--Vinogradov Theorem cited
earlier.

\proclaim{Theorem 4} Theorem 3 holds unconditionally if the field
$k$ is contained in a Galois\break extension $K$ such that either:
 \vskip 6pt
 {\bf (a)} ${\frak{g}} = \text{\rm
Gal}(K/\Bbb Q)$ has an abelian subgroup of index at most 4 (in
particular, if $k$ is an abelian extension),

or
 
 {\bf (b)}
$[K:\Bbb Q] < 42.$
\endproclaim

The proof of the upper bound is very different from the proof of
the lower bound. For a group $A$, we denote by $s_r(A)$, the
number of subgroups of $A$ of index at most $n$. A somewhat
involved reduction process is applied to show that the problem of finding the upper bound is
actually equivalent to an extremal counting problem of subgroups
of finite abelian groups (see section \S 5) which is given in Theorem 5. A sharp upper bound for that counting problem follows from
the case $R = 1$ of the following theorem. 

\proclaim{Theorem 5} Let $R \ge 1$ be a real number and let d be a fixed integer $\ge 1.$
Suppose that $A = C_{x_1}\times C_{x_2}\times \cdots\times C_{x_t}$ is
an abelian group such that the orders $x_1, x_2, \ldots, x_t$ of
its cyclic factors do not repeat more than $d$ times each. Suppose
that $r|A|^R \le n$ for some positive integers $r$ and $n$. Then
as $n$ tends to infinity, we have
$$s_r(A) \le n^{(\gamma+o(1))\ell(n)},$$
where $\gamma = \frac{(\sqrt{R(R+1)}-R)^2}{4R^2}.$
\endproclaim

In an earlier version of this paper, Theorem 5 was proved in a similar manner, but only for $R = 1.$ The more general case was proved in an early version of \cite{LN}. We thank the authors of \cite{LN} for allowing us to include the general version here.

The above results suggest that for every Chevalley group scheme
$G$, the upper and lower limiting constants, $\alpha_\pm (G(\Cal
O_S))$ are equal to each other, and depend only on $G$ and not on
$\Cal O$. In fact, we can make a precise conjecture, for which we
need to introduce some additional notation. Let $G$ be a Chevalley
group scheme of dimension $d=\dim(G)$ and rank $\ell=rk(G)$. Let
 $\kappa = \vert\Phi^+\vert$ denote the number of positive roots in
the root system of $G$, and  let $R = R(G) =
\frac{d-\ell}{2\ell}=\frac{\kappa}{\ell}$. We see that if $G$ is of type
$A_\ell$ (resp. $B_\ell, C_\ell, D_\ell,
G_2,F_4,E_6,E_7,E_8$)
then $R =
\frac{\ell +1}{2}, \quad(\text{resp.} \ell, \ell, \ell-1, 3, 6,
6, 9, 15)$.

\proclaim{Conjecture}  Let $k, \Cal O$, and $S$ be as in Theorem 3,
 and suppose  that $G$ is a simple Chevalley group scheme.  Then
$$
\alpha_+ (G(\Cal O_S)) = \alpha_- (G(\Cal O_S)) = \frac{\left(
\sqrt{R(R+1)}-R\right)^2}{4R^2}.
$$
\endproclaim

The conjecture reflects the belief that ``most'' subgroups of
$H=G(\Bbb Z/m\Bbb Z)$ lie between the Borel subgroup $B$ of $H$
and the unipotent radical of $B$. We prove here the lower bound of
the general conjecture (under the same assumptions as in Theorem 3
and 4). In our earlier version this was done only for Galois
extensions, but it was observed in an earlier version of \cite{LN}
that a small modification of the argument works in the general
case. We thank the authors of \cite{LN} for allowing us to make these small modifications
 here.

This paper gives a complete proof of the upper bound for the case
of $SL_2$, based on the known detailed classification of subgroups
of $SL_2(\Bbb F_q)$ for finite fields $\Bbb F_q$ of order $q$. We
also give a partial result towards the upper bound in the general
case. The upper bound is proved in full for every field $k$ in
\cite{LN}. The reader is also referred to a more general version
there when $G$ is not assumed to be split.

\proclaim{Theorem 6}  Let $k, \Cal O,$ and $S$ be as in Theorem 3. Let
$G$ be a simple Chevalley group scheme of dimension $d$ and rank
$\ell$, and $R=R(G)=\frac{d-\ell}{2\ell}$,   then: \vskip 10pt
{\bf (a)} Assuming $GRH$ or the assumptions of Theorem 4,
$$\alpha_- (G(\Cal O_S)) \geq
\frac{\left(\sqrt{R(R+1)}-R\right)^2}{4R^2}\sim \frac{1}{16 R^2}.$$
\vskip 1pt {\bf (b)} There exists an absolute constant $C$ such
that $$\alpha_+ (G(\Cal O_S))\leq \; C\cdot \frac{
\big(\sqrt{R(R+1)}-R\big)^2}{4R^2}.$$
\endproclaim

\noindent {\bf Remark:} As the upper bound is proved in full in
\cite{LN} (i.e., $C = 1$ in part (b)) we omit in this paper the proof of part (b) of Theorem
6.

\proclaim{Corollary 7} There exists an absolute constant $C$ such that
for $d = 2, 3, \ldots$
$$
(1-o(1))\frac{1}{4d^2} \; \leq \alpha_- (SL_d(\Bbb Z)) \leq
\alpha_+ (SL_d(\Bbb Z))\leq C \frac{1}{d^{2}}.
$$
\endproclaim
This greatly improves the upper bound $\alpha_+ (SL_d(\Bbb Z)) < \frac{5}{4}
d^2$ implicit in \cite{Lu}
and settles a question asked there.

As a byproduct of the proof of Theorem 5 in \S6 we obtain the
following. \proclaim{Corollary 8} The subgroup growth type of
$SL_d(\Bbb Z _p)$  is at least $n^c$ where $$c=(3-2\sqrt{2})d^2-
2(2-\sqrt{2}),$$
and where $\Bbb Z _p$ denotes the ring of $p$-adic integers.
\endproclaim

  The counting techniques in this paper can be applied to solve a novel extremal
problem in multiplicative number theory involving the greatest
common divisors of pairs $(p-1,  p'-1)$ where $p, p'$ are prime
numbers. The solution of this problem does not appear amenable to
the standard techniques used in  analytic number theory. Considering this problem first was crucial for obtaining Theorem 5.

\proclaim{Theorem 9} For $n \to \infty$, let
$$M(n) = \max \Bigg\{ \prod_{p, p' \in \Cal P} \text{\rm gcd}(p-1, \; p'-1) \;\, \Bigg | \;\, \Cal P = \text{set of distinct primes
where} \, \prod_{p\in\Cal P} p \le n\Bigg\}.$$ Then we have: $$\lim\limits_{n\to\infty} \; \frac{\log M(n)}{\lambda(n)} \; = \; \frac14,$$  where
$\lambda(n) = (\log n)^2/\log\log n)$.
\endproclaim

\vskip 10pt
The paper is organized as follows.

In \S1, we present some required preliminaries and notation.

In \S2, we introduce the notion of a Bombieri set which is the
crucial ingredient needed in the proof of the lower bounds. We
then use it in \S3 and \S4 to prove the lower bounds of Theorems
2, 3,  4, and 6. We then turn to the proof of the upper bounds. In
\S5, we show how the counting problem of congruence subgroups in
$SL_2(\Bbb Z)$ can be completely reduced to an extremal counting
problem of subgroups of finite abelian groups;  the problem is
actually, as one may expect, a number theoretic extremal problem -
see \S6 and \S7 where this extremal problem is solved and the
upper bounds of Theorems 2, 3, and 4 are then deduced in \S8.
Finally, in \S9 we prove Theorem 9.

  The results of this paper are announced in \cite{GLNP}.

\pagebreak

\vskip .20in
{\bf \S1.  Preliminaries and notation}
\vskip .10in

Throughout this paper we let $$\ell (n) = \frac{\log n}{\log\log
n}, \qquad\lambda(n) = \frac{(\log n)^2}{\log\log n}.$$ All logarithms in this paper are to base $e$.
  If $f$ and $g$ are
functions of $n$, we will say that {\it $f$ is small w.r.t. $g$}
if $\lim\limits_{n\to\infty} \frac{\log f(n)}{\log g(n)} = 0$. We
say that $f$ is {\it small} if $f$ is {\it small} with respect to
$n^{\ell (n)}$.  Note that if $f$ is small, then multiplying
$C_n(\Gamma)$ by $f$ will have no effect  on the estimates of
$\alpha_+(\Gamma)$ or $\alpha_- (\Gamma)$.  We may, and we will,
ignore factors which are small.

Note also that if $\varepsilon (n)$ is a function of $n$ which
is smaller than $n$ \newline (i.e., $\log \varepsilon (n) =
o(\log n)$) then:
$$\overline{\lim} \; \frac{\log C_{n\varepsilon
(n)}(\Gamma)}{\lambda (n)} = \alpha_+ (\Gamma) \tag1.1
$$
and
$$\underline{\lim} \; \frac{\log C_{n\varepsilon
(n)}(\Gamma)}{\lambda (n)} = \alpha_- (\Gamma). \tag1.2
$$

The proof of (1.1) follows immediately from the inequalities:
 $$\align\alpha_{+}(\Gamma) &= \;\overline{\lim} \; \frac{\log C_{n}(\Gamma)}{\lambda (n)}
\; \le \; \overline{\lim} \; \frac{\log
C_{n\varepsilon (n)}(\Gamma)}{\lambda(n)}\\
 &= \;\overline{\lim} \; \frac{\log
C_{n\varepsilon (n)}(\Gamma)}{\lambda (n\varepsilon
(n))}\cdot \frac{\lambda (n\varepsilon (n))}{\lambda(n)} \\
&\leq \;
\alpha_+(\Gamma )\cdot 1\\
 &= \;\alpha_+ (\Gamma).
\endalign
$$
Here, we have used the fact that $\overline{\lim} \;
\frac{\lambda (n\varepsilon (n))}{\lambda (n)} = 1$, which is
an immediate consequence of the assumption that
$\varepsilon (n)$ is small with respect to $n$. A similar argument proves
(1.2).

It follows that we can, and we will sometimes indeed, enlarge
$n$ a bit when evaluating $C_n(\Gamma)$, again without
influencing $\alpha_+$ or $\alpha_-$.  Similar remarks
apply if we divide $n$ by $\varepsilon (n)$ provided
$\varepsilon (n)$ is bounded away from $0$.

The following lemma is proved in \cite{Lu} in a slightly weaker form
and in  its current form is proved in \cite{LS, Proposition 5.1.1}.

\proclaim{Lemma 1.1} (``Level versus index'').
   Let $\Gamma$ be as before.  Then there exists a constant $c > 0$
such that if $H$ is a congruence subgroup of $\Gamma$ of index at
most $n$, then $H$ contains $\Gamma (m)$ for some $m\leq cn$, where
$m\in \Bbb Z$ and by $\Gamma(m)$ we mean $\Gamma (m\Cal O_S)$.
\endproclaim

\proclaim{Corollary 1.2}  Let $\gamma_n(\Gamma) =
\sum\limits^n_{m=1} s_n(G(\Cal O_S/m\Cal O_S))$, where for a group
$H$, $s_n(H)$ denotes the number of subgroups of $H$ of index at
most $n$.  Then we have $\alpha_+ (\Gamma) = \overline{\lim} \;
\frac{\log \gamma_n(\Gamma)}{\lambda (n)}$ and $\alpha_- (\Gamma)
= \underline{\lim} \; \frac{\log \gamma_n(\Gamma)}{\lambda (n)}$.
\endproclaim

\demo{Proof}  By Lemma 1.1,  $C_n(\Gamma ) \leq
\gamma_{cn} (\Gamma)$ for some $c > 0$.  It is also clear that
 $\gamma_n(\Gamma) \leq n\cdot C_n(\Gamma)$.
 Since  $c$ is small w.r.t. $n$, Corollary 1.2
follows by
arguments of the type we have given above.  \qed
\enddemo

The number of elements in a finite set $X$ is denoted by $\#X$ or
$|X|$. The set of subgroups of a group $G$ is denoted by $\text{Sub}(G)$.

\vskip .20in
{\bf \S2.} {\bf Bombieri Sets.}
\vskip .10in

We introduce some  additional notation. Let $a,q$ be relatively prime integers with $q > 0$.  For
$x > 0$, let $\Cal P(x;q,a)$ be the set of primes $p$ with $p\leq
x$ and $p\equiv a  \pmod q$.  For $a=1$, we set $\Cal P(x;q) =
\Cal P(x;q,1)$.  We also define $$\vartheta (x;q,a) =
\sum\limits_{p\in\Cal P(x;q,a)} \log p.$$

If $f(x), \, g(x)$ are arbitrary functions of a real variable
$x$, we say
$f(x) \sim g(x)$ as $x \to \infty$ if
$$
\lim_{x\to\infty} \,
\frac{f(x)}{g(x)} \; = \; 1.
$$

Define the error term
$$
E(x; q, a) = \vartheta(x; q, a) - \frac{x}{\phi(q)},
$$
where $\phi(q)$ is Euler's function.
Then Bombieri proved the following deep theorem
\cite{Bo}, \cite{Da}.

\proclaim{Theorem 2.1} (Bombieri) Let $A > 0$ be
fixed. Then
there exists a constant $c(A) > 0$ such that
$$\sum_{q \;\le\; \frac{\sqrt{x}}{(\log x)^{A}} }\max_{y \le
x}\max_{(a, q) = 1} \big |E(y; q, a)\big| \; \le \;
c(A)\cdot\frac{x}{(\log x)^{A-5}}$$
as $x \to \infty.$
\endproclaim

 This theorem shows that the error terms
$\underset {(a, q) =
1}\to{\max}\; E(x; q, a)$ behave as if they satisfy the
Riemann hypothesis in an averaged sense.
\vskip .10in

\proclaim{Definition 2.2}  Let $x$ be a large positive real
number.  A {\bf Bombieri prime} (relative to $x$) is a prime
$q \le \sqrt{x}$ such that the set
$\Cal P(x, q)$ of primes $p \le x$ with
$p \equiv 1 \pmod{q}$ satisfies
$$\max_{y \;\le \; x} \; |E(y; q, 1)|  \;\le \;
\frac{x}{\phi(q) (\log x)^2}.$$ We call $\Cal P(x, q)$ a
{\bf Bombieri set} (relative to $x$).
\endproclaim

\noindent {\bf Remark}. In all the applications in this paper, we
do not really need $q$ to be prime, though it makes the
calculations somewhat easier.  We could work with $q$ being a
``Bombieri number".

\proclaim{Lemma 2.3} Fix $0 < \rho < \frac12.$ Then for $x$
sufficiently large, there exists at least one Bombieri prime
(relative to $x$)  $q$ in the interval
$$\frac{x^\rho}{\log x} \le q \le x^\rho.$$
\endproclaim

\demo{Proof} Assume that
$$
\max_{y \;\le \; x} |E(y; q, 1)|  \; >\;
\frac{x}{\phi(q) (\log x)^2}
$$
for all primes $\frac{x^\rho}{\log x} \le q \le x^\rho,$ i.e.,
that there are no such Bombieri primes  in the interval. In view
of the trivial inequality, $\phi(q) = q - 1 < q,\,$ it immediately
follows that
$$\sum_{\frac{x^\rho}{\log x} \,\le \, q \,\le \, x^\rho}
\max_{y \;\le \; x}\big |E(y; q, 1)\big| \; > \; \frac{x}{(\log x)^2}
\sum_{\frac{x^\rho}{\log x} \,\le \, q \,\le \, x^\rho} \frac{1}{q} \; > \;
\frac{x\cdot (\log\log x)}{2\rho\cdot (\log x)^3},
$$
say, for sufficiently large $x$.  This follows from the well
known asymptotic formula \cite{Lan} for the partial sum of the
reciprocal of the primes
$$\sum_{q \,\le \, Y} \; \frac{1}{q} \; = \;  \log\log Y + b +
O\left(\frac{1}{\log Y}\right)$$ as
$Y \to \infty.$ Here $b$ is an absolute constant.  This
contradicts
Theorem 2.1  with
$A \ge 8$ provided $x$ is sufficiently large. \qed
\enddemo
\vskip .10in

\proclaim{Lemma 2.4}  Let $\Cal P(x, q)$ be a Bombieri set.
Then for $x$ sufficiently large
$$
\left |\#\Cal P(x, q)\; - \; \frac{x}{\phi(q) \log x} \right | \;
\le \;  3\left(\frac{x}{\phi(q) (\log x)^2}\right).
$$
\endproclaim

\demo{Proof} We have
$$
\align \sum_{p\in \Cal P(x, q)} 1 \; &= \; \sum_{n=2}^x
\frac{\vartheta(n; q, 1) - \vartheta (n-1; q , 1)}{ \log n}\\
&= \sum_{n=2}^x \vartheta(n; q, 1) \Big(\frac{1}{\log(n)} -
\frac{1}{\log(n+1)}\Big )
 +
\frac{\vartheta(x; q, 1)}{ \log ([x]+1)}\\
& = \sum_{n=2}^x \vartheta(n; q, 1)
\frac{\log\left(1+\frac1n\right)}{\log n\log(n+1)} +
\frac{\vartheta(x; q, 1)}{\log x} - \vartheta(x; q, 1)
\left(\frac{1}{\log x} - \frac{1}{\log([x]+1)}\right).\endalign
$$

It easily follows that
$$\left |\sum_{p\in \Cal P(x, q)} 1 -
\frac{\vartheta(x; q, 1)}{\log
x}\right | \; \le \; \sum_{n=2}^x \vartheta (n; q, 1)
\frac{1}{n\cdot(\log n)^2} + \vartheta(x; q, 1)
\left (\frac{1}{\log x} -
\frac{1}{\log(x+1)}\right ).
$$
By the property of a Bombieri set, we have
the estimate $|\vartheta(n; q, 1) -
\frac{n}{\phi(q)}| \le \frac{x}{\phi(q) (\log x)^2}$, for $n \le x$.
Since $\left (\frac{1}{\log x} -
\frac{1}{\log(x+1)}\right ) = \frac{\log\left(1 + \frac{1}{x}\right)}{\log x\log (x+1}) =
O\left( \frac{1}{x(\log x)^2} \right)$, the second expression on the right side of the above equation is very small and
can be ignored. It remains to estimate the sum $\sum\limits_{n=2}^x \vartheta (n; q,
1)
\frac{1}{n\cdot(\log n)^2}$. This sum can be broken into two parts, the first of which
corresponds to $n \le \frac{x}{(\log x)^3},$  which is easily seen to be very small, so
can be ignored. We estimate
$$\align\sum_{\frac{x}{(\log x)^3} \le n\le x} \vartheta (n; q, 1)
\frac{1}{n\cdot(\log n)^2} \; &\; = \;  \sum_{\frac{x}{(\log x)^3} \le n\le x}
\frac{n}{\phi(q)}
\cdot
\frac{1}{n(\log n)^2}\\
&\hskip 65pt \; + \;  O\left(\sum_{\frac{x}{(\log x)^3} \le n\le x}
\frac{x}{\phi(q)(\log x)^2}
\cdot
\frac{1}{n(\log n)^2}\right)\\ & = \;  \sum_{\frac{x}{(\log x)^3} \le n\le x}
\frac{1}{\phi(q) (\log n)^2} \; +  \;  O\left(\frac{x}{\phi(q) (\log x)^3}  \right) \\
& \le \;
\frac{3}{2} \frac{x}{\phi(q)(\log x)^2},\endalign$$
which holds for $x$ sufficiently large and where  the constant
$\frac32$ is not optimal. Hence
$$\left |\sum_{p\in \Cal P(x, q)} 1 -
\frac{\vartheta(x; q, 1)}{\log
x}\right | \; \le \; \frac{7}{4} \frac{x}{\phi(q)(\log x)^2},$$
say.
 Since $|\vartheta(x; q, 1) - \frac{x}{\phi(q)}| \le
\frac{x}{\phi(q) (\log x)^2},$ Lemma 2.4 immediately follows.\qed
\enddemo

\vskip 20pt
{\bf \S3.} {\bf Proof of the lower bound over $\Bbb Q$.}
 \vskip .10in
In this section we consider the case of $k=\Bbb Q$ and
$\Cal O = \Bbb Z$.

Fix a real number $0 < \rho_0 < \frac12.$ It follows from Lemma 2.3
that for  $x \to \infty$ there exists a real number $\rho$ which converges to
$\rho_0$, and a prime number
$q
\sim x^\rho$ such that $\Cal P(x, q)$ is a Bombieri set.

Define
$$
P = \prod_{p \; \in \;\Cal P(x,q)} p.
$$
It is clear from the definition of a Bombieri set that
$$
\log P \; \sim \; \frac{x}{\phi(q)} \sim x^{1-\rho}
$$
and from Lemma 2.4 that
$$ L = \#\Cal P (x, q) \sim \frac{x}{\phi(q)\log x}
\sim\frac{x^{1-\rho}}{\log x}.
$$

Consider $\Gamma(P)= \ker (G(\Bbb Z)\to G(\Bbb Z / P\Bbb Z)) $
which is of index at most $P^{\text{dim}(G)}$ in $\Gamma.$ Note
that for every subgroup $H/\Gamma(P)$ in $\Gamma/\Gamma(P)$ there
corresponds a subgroup $H$ in $\Gamma$ of index at most
$P^{\text{dim}(G)}$ in $\Gamma$.

 By strong approximation
$$\Gamma/\Gamma(P) = G\left(\Bbb Z/P\Bbb Z\right ) \cong
\prod_{p\in\Cal P(x, q)} G(\Bbb F_p).$$

  Let $B(p)$ denote the Borel subgroup  in $G(\Bbb F_p).$ Then
$$\log\Big(\#B(p)\Big) \; \sim \;
\frac{\text{dim}(G) +
\text{rk}(G)}{2}\, \log p,$$
where rk$(G)$ denotes the rank of $G$ as an algebraic group.
But
$$\log\Big(\#G(\Bbb F_p)\Big)\; \sim \;
\text{dim}(G)\, \log p .$$
It immediately follows that (for $p \to \infty$)
$$\log\big [G(\Bbb F_p) : B(p)\big] \; \sim \;
\frac{\text{dim}(G) - \text{rk}(G)}{2}\, \log p,$$
 and, therefore,
$$
\log\big [G(\Bbb Z/P\Bbb Z) : B(P)\big ] \sim
\frac{\text{dim}(G) - \text{rk}(G)}{2}\, \log P.
$$
where $B(P)\le G(\Bbb Z/P\Bbb Z)$ is:
$$B(P) = \prod\limits_{p\in\Cal P(x:q)} B(\Bbb F_p).
$$

Now $B(p)$ is mapped onto ${\Bbb F_p^\times}^{\text{rk}(G)}$ and,
hence, is also mapped onto $\left(\Bbb Z/q\Bbb
Z\right)^{\text{rk}(G)}$ since $\#F_p^\times = p-1$ and $p \equiv
1 \pmod{q}.$ So $B(P)$ is mapped onto
$$\left(\Bbb Z/q\Bbb Z\right)^{\text{rk}(G)\cdot L}$$
where
$$L = \#\Cal P(x, q) \sim \frac{x}{\phi(q)\log x}
\sim \frac{x^{1-\rho}}{\log x}.$$

 For a real number $\theta$, define $\lceil \theta\rceil$ to be
the smallest integer $t$ such that $\theta \le t.$ Let
$0 \le \sigma \le 1.$

We will now use Proposition 6.1, a basic result on counting
subspaces of finite vector spaces. It follows
that $B(P)$ has at least
$$q^{\sigma (1-\sigma) \text{rk}(G)^2 L^2 + O(\text{rk}(G)\cdot L)}$$
subgroups of index equal to
$$\iota = 
q^{\lceil \sigma\cdot \text{rk}(G)\cdot L\rceil} \cdot
\big [G(\Bbb Z/ P\Bbb Z) :
B(P)\big].
$$
Hence, for $x \to \infty$,
$$\align
&\log\Big(\text{\#\big\{subgroups\big\}} \Big) = \Big(\sigma
(1-\sigma) \text{rk}(G)^2 L^2 + O(\text{rk}(G)\cdot L)\Big)
\log q  \\
&\qquad \sim \;\; \sigma (1-\sigma) \text{rk}(G)^2
\frac{x^{2-2\rho}}{(\log x)^2} \cdot\rho\log x,
\endalign
$$
while
$$
\align \log(\iota) &= \lceil \sigma\cdot\text{rk}(G)
\cdot  L\rceil\cdot\log q  +\frac12 \big(\text{dim}(G) -
\text{rk}(G)\big)\log P\\ &\sim \text{rk}(G)\sigma
\frac{x^{1-\rho}}{\log x} \rho \log x \, + \,
\frac12 \big(\text{dim}(G) - \text{rk}(G)\big) x^{1-\rho}\\
&= \Big(\sigma\cdot\rho \cdot\text{rk}(G) + \frac12
\big(\text{dim}(G) -
\text{rk}(G)\big)\Big) x^{1-\rho},
\endalign
$$
and
$$
\log\log(\iota) \sim (1-\rho) \log x.
$$
It is clear from the estimate for $\log \iota$ above that given any
index $n>>0$ we can choose $x$ such that $\log \iota \sim \log n$.
We  compute
$$\align\frac{\log\Big(\text{\#\{subgroups\}}
\Big)}{\big(\log(\text{index})\big)^2/\log\log(\text{index})}
\;   &\sim   \;
\;\frac{ \sigma (1-\sigma)\cdot\text{rk}(G)^2\cdot\rho \;
\frac{x^{2-2\rho}}{\log
x}   } {\Big(\Big(\sigma\cdot\rho \cdot \text{rk}(G) +
\frac12 \big(\text{dim}(G) -
\text{rk}(G)\big)\Big) x^{1-\rho}\Big)^2\Big / (1-\rho) \log x} \\
& \\
& \sim \;  \frac{\sigma  (1-\sigma)  \rho  (1-\rho) \cdot
\text{rk}(G)^2}{\Bigl(\left(\sigma \rho - \frac12\right)\cdot
 \text{rk}(G)  +
\frac12 \,\text{dim}(G)\Bigr)^2}\endalign$$

as $x \to \infty.$

  We may rewrite
$$\frac{\sigma  (1-\sigma)  \rho  (1-\rho) \cdot
\text{rk}(G)^2}{\Bigl(\left(\sigma \rho - \frac12\right)\cdot
 \text{rk}(G)  +
\frac12 \,\text{dim}(G)\Bigr)^2} \;\, = \;\, \frac{\sigma(1-\sigma)
\rho(1-\rho)}{(\sigma\rho + R)^2}$$
where
$$R = \frac{\dim(G) - \text{rk}(G)}{2\cdot\text{rk}(G)}.$$

Now, for fixed $R$, it is enough to choose
$\sigma, \rho$ so
that
$$ \frac{\sigma(1-\sigma) \rho(1-\rho)}{(\sigma\rho + R)^2}$$
is  maximized. This occurs when
$$\rho = \sigma = \sqrt{R(R+1)} - R,$$
in which case we get
$$ \frac{\sigma(1-\sigma) \rho(1-\rho)}{(\sigma\rho + R)^2} \; =
\frac{\left(\sqrt{R(R+1)}-R\right)^2}{4R^2}.$$

Actually, we choose $\rho_0$ to be $\sqrt{R(R+1)}-R$, then we can
take $\rho$ to be asymptotic to $\rho_0$ as $x$ is going to
infinity. Note that
$\frac{\left(\sqrt{R(R+1)}-R\right)^2}{4R^2}<\frac{1}{16R^2}$
holds for all $R>0$. This follows from the easy inequality
$\sqrt{R(R+1)}-R\le \frac12$. It is also straightforward to see
that  $\sqrt{R(R+1)}-R$ converges to $\frac12$ as $R\to\infty$
hence
$\frac{\left(\sqrt{R(R+1)}-R\right)^2}{4R^2}\sim\frac{1}{16R^2}$.

In the special case
when $R = 1$, we obtain the lower bound of Theorem 2. For a
simple Chevalley group   scheme over $\Bbb Q$, this gives the
lower bound in Theorem 6.

\vskip 20pt
{\bf \S 4. Proof of the lower bound for a general number field.}
\vskip .10in

To prove the lower bounds over a general number field we need an
extension of the Bombieri--Vinogradov Theorem to these fields, as was obtained by
Murty and Murty \cite{MM}.

Let us first fix some notations:

Let $k$ be a finite  extension of degree $f$ over $\Bbb Q$, $K$
its Galois closure of degree $d$, ${\frak{g}} = \text{Gal}(K/\Bbb
Q),$ and $\Cal O_k$ the ring of integers in $k$. For a rational
prime $q $ and $x \in\Bbb R$, we will denote by $\tilde \Cal
P_K(x, q)$ the set of rational primes $p\equiv 1 (\mod q)$ where
$p$ splits completely in $K$ and $p \le x$.  Let
$$\tilde \pi_K(x, q)= \#
\tilde\Cal P_K(x, q), \qquad \tilde \nu_K(x, q) =
\sum\limits_{p\in \tilde\Cal P_K(x, q)} \log p,$$ and,
$$
\tilde E_K(x, q)= \tilde \nu_K(x, q) - \frac{x}{d\phi (q)}.
$$
We shall show that the following theorems follow from Murty and Murty \cite{MM}.

\proclaim{Theorem 4.1} Let $K$ be a fixed finite Galois extension
of $\Bbb Q$. Assume GRH (generalized Riemann hypothesis) for $K$
and all cyclotomic extensions $K(\zeta_\ell)$ with $\ell$ a
rational prime and $\zeta_\ell$ a primitive $\ell^{\text{\rm th}}$
root of unity. Then for every $0 < \rho < \frac{1}{2}$ and $x \to
\infty$, there exists a rational prime $q$ such that \vskip 8pt
{\bf (a)} $\frac{x^\rho}{\log x} \le q \le x^\rho$

\vskip 10pt {\bf (b)} $|\tilde\pi_K(x, q)-\frac{x}{d'\phi(q)\log x}
|\le 3\left(\frac{x}{d'\phi(q)(\log x)^2}\right)$

\vskip 10pt  {\bf (c)} $\mathop{\max}\limits_{y\le x} |\tilde
E_K(y, q) |\le \frac{x}{d'\phi(q)(\log  x)^2},$
\vskip 8pt\noindent
where $d' = [K : Q]/t$ and $t$ denotes the degree of the intersection of $K$ and
the cyclotomic field $\Bbb Q(\zeta_q)$ over $\Bbb Q$.

\endproclaim
\vskip 8pt
\noindent
{\bf Remark:} In fact, GRH gives a stronger result than what is stated in Theorem 4.1. For example, it can be shown 
 that for every prime $q < x^{\frac12}$ the error terms in parts (b), (c), take the form $\Cal O\left (x^{\frac12} \log(qx)\right)$
 (see \cite{MMS} for a more precise bound). Theorem 4.1 is stated in this special form because it can be proved unconditionally in some cases.

\vskip 15pt \proclaim{Theorem 4.2} Theorem 4.1 can be proved
unconditionally for $K$ if either: \vskip 6pt {\bf (a)} ${\frak{g}}
= \text{\rm Gal}(K/\Bbb Q)$ has an abelian subgroup of index at
most 4 (this is true, for example, if $k$ is an abelian
extension); \vskip 3pt \noindent
or
\vskip 3pt
{\bf (b)} $[K:\Bbb Q] <
42.$
\endproclaim

\vskip 10pt \proclaim{Theorem 4.3} Theorem 4.1 is valid
unconditionally for every $K$ with the additional\break assumption
that $0 < \rho < \frac{1}{\eta}$, where $\eta$ is the maximum of 2
and $d^*- 2$, and where $d^*$ is  the index of the largest
possible abelian subgroup of ${\frak{g}} = \text{\rm Gal}(K/\Bbb
Q).$ In particular, we may take $\eta = d^* - 2$ if $d^* \ge 4 $
and $\eta = 2$ if $d^* \le 4$.
\endproclaim

\vskip 10pt \noindent{\it Proof of Theorems 4.1 - 4.3.} For any $\epsilon >
0, A > 0,$ under the
assumptions of Theorem 4.1 or 4.2 (a), Murty and Murty \cite{MM} prove the
following Bombieri theorem:
$$\sum_{q \le x^{\frac12-\epsilon}} \max_{(a,q)=1} \max_{y\le x} \left |\pi_C(y,q,a) -
\frac{|C|}{|G|}\cdot\frac{1}{\phi(q)}
\pi(y)\right | \ll \frac{x}{(\log x)^A}.\tag 4.1
$$
Here $C$ denotes a conjugacy class in $\frak g$,
$\pi(y) = \sum_{p \le y} 1,$
$$\pi_C(x,q,a) = \underset p\; \text{unramified in $K$}\to{\underset p \equiv a \pmod{q}\to
{\underset (p, K/\Bbb Q) = C \to{ \sum_{p \le x}}}} 1,$$  and $(p,
K/\Bbb Q)$ denotes the Artin symbol.

In fact, under the assumption of the $GRH$, equation (4.1) holds, but without
assuming
$GRH$ they showed that (4.1) holds when the sum is over
$q < x^{\frac{1}{\eta} - \varepsilon}$ where $\eta$ is defined as
follows: Let
$$
d^* = \min_{H} \max_{w} [\frak g:H] w(1) \tag 4.2
$$

The minimum here is over all subgroups $H$ of $\text{Gal}(K/Q)$
satisfying:
\vskip 6pt
{\bf (i)} $H\cap C\neq \emptyset$, and
\vskip 6pt
{\bf (ii)} for every irreducible character $w$ of $H$ and any
non-trivial Dirichlet character $\chi$, the Artin $L$-series $L(s,
w\otimes \chi)$ is entire. 
\vskip 6pt
\noindent
 Then the  maximum in (4.2) is over the
irreducible characters of such $H$'s.

Now
$$
\eta = \cases\matrix d^* - 2 & \hbox{\rm if} \; \;\;  d^* \ge 4\\
2 &\hbox{\rm if} \; \;\; d^* \le 4\endmatrix\endcases
$$

We need their result for the special case when $C$ is the identity
conjugacy class.  In this case $\frac{|C|}{| \frak g |} = \frac {1}{d'}$ and
$\pi_C (y, q, 1) = \tilde \pi_k (y, q)$.  So for proving Theorem
4.3 we can take for $H$ an abelian subgroup of smallest index and
then $H$ satisfies assumption (i) and (ii).  (Recall that abelian
groups satisfy (AC) - Artin conjecture, i.e. $L(s, w \otimes
\chi)$ are entire -- see \cite{CF}).

For Theorem 4.2(a), again take $H$ to be the abelian subgroup of
index at most 4.  It satisfies (i) and (ii) and this time $\eta =
2$.

For Theorem 4.2(b): Going case by case over all possible numbers
$d < 42$, one can deduce by elementary group theoretic arguments
that every finite group $\frak{g}$ of order $d < 42$, has an
abelian subgroup of index at most 4, unless $d = 24$ and $\frak g$
is isomorphic  to the symmetric group $S_4$.   But for this group,
Artin \cite{CF} proved Artin's conjecture in 1925.  Moreover, every irreducible
character of $S_4$ is of degree at most 4.  Thus for
${\frak{g}}=S_4$ we have $d^* = 4$ and so $\eta = 2$.

The proofs of Theorems 4.1, 4.2 and 4.3 follow now in the
same manner as in \S 2.

  Using Theorems 4.1, 4.2, 4.3, we
can now prove the lower bounds of Theorem 3 and 4 just as in \S 3.
Note that for every prime $p \in \tilde\Cal P_K (x, q)$ we may
take an ideal  $\pi = \pi(p)$ in $\Cal O_k$ with $[\Cal O_k:
\pi]=p, \pi\cap\Bbb Z = p\Bbb Z$.  Let
$$P = \prod\limits_{p\in \tilde\Cal P_K(x,q)} \pi(p).$$ 
Then, since  $x \to \infty$, we may choose $q, \rho$ (using Theorem 4.1)  so that
 $$\log [\Cal O
:P]\sim\frac{x}{d\, \phi(q)} \sim \frac{x^{1-\rho}}{d}, \qquad L:=|
\Cal P_K(x, q)|\sim\frac{x}{d\,\phi(q)\log x} \sim \frac{x^{1-\rho}}{d \log x},$$ and
$$G(\Cal O/P) = \prod\limits_{p\in \Cal P_K(x,q)} G(\Cal O/\pi(p))
\simeq \prod\limits_{p\in \tilde\Cal P_k(x, q)} G(\Bbb Z/p \Bbb
Z).$$ We can now take for every rational prime $p \in\tilde \Cal
P_k(x, q)$, the Borel subgroup $B(p)$ as in \S 3 and define:
$$
B(P)=\prod\limits_{p\in\tilde\Cal P_k(x, q)} B(p).
$$
Then $B(P)$ is mapped onto $(\Bbb Z/q\Bbb Z)^{rk(G)\cdot L}$ and
$$\log\big [G(\Cal O/P): B(P)\big ]\sim
\frac{\dim(G)-rk(G)}{2}\cdot\log [\Cal O
:P].$$ Thus, by a computation similar to the one  in \S 3 (note that the $d$'s cancel in this computation), we can show that
$$
\alpha_-(G(\Cal O))\ge \frac{\left(
\sqrt{R(R+1)}-R\right)^2}{4R^2}.
$$
The lower bounds of Theorems 3, 4, and 6 are now also proved. We
now turn to the proof of the upper bounds.

\vskip .2in
{\bf \S5.  From $SL_2$ to abelian groups}
\vskip .10in

In this section we show how to reduce the estimation of
$\alpha_+(SL_2(\Bbb Z))$ to a problem on abelian groups.

Corollary 1.2 shows us that in order to give an upper bound on
$\alpha_+(\Gamma)$ it suffices to bound $s_n(G(\Bbb Z/m\Bbb Z))$
when $m\leq n$.  Our first  goal is to show that we can further
assume that $m$ is a product of different primes. To this end
denote $\overline m = \prod p$ where $p$ runs through all the
primes dividing $m$.

We have an exact sequence
$$
1\to K\to G(\Bbb Z/m\Bbb Z)\overset{\pi}\to \longrightarrow G(\Bbb
Z/\overline{m}\Bbb Z)\to 1
$$
where $K$ is a nilpotent group of rank at most $\dim G$.
Here, the rank of a finite group $G$ is defined to be the smallest integer $r$ such that
every subgroup of $G$
is generated by $r$ elements, (see \cite{LS, Window 5, \S 2}).

\proclaim{Lemma 5.1}  Let $1\to K\to U\overset{\pi}\to
\longrightarrow L\to 1$ be an exact sequence of finite groups,
where $K$ is a solvable group of derived length $\ell$ and of
rank at most $r$.  Then the number of supplements to $K$ in $U$
(i.e., of subgroups $H$ of $U$ for which $\pi (H)= L$) is bounded
by $\vert U\vert^{3r^2+\ell r}$.
\endproclaim

\demo{Proof} See \cite{LS, Corollary 1.3.5}.
\enddemo

 \proclaim{Corollary 5.2}  $s_n(G(\Bbb Z/m\Bbb Z)) \leq m^{f'(\dim
G)\log\log  m}\cdot s_n(G(\Bbb Z/\overline m\Bbb Z))$ where $f'(\dim
G)$ depends only on $\dim G$.
\endproclaim

\demo{Proof}  Let $H$ be a subgroup of index at most $n$ in
$G(\Bbb Z/m\Bbb Z)$ and denote $L = \pi (H)\leq G(\Bbb
Z/\overline m\Bbb Z)$.  So $L$ is of index at most $n$ in $G(\Bbb
Z/\overline m\Bbb Z)$.  Let $U = \pi^{-1}(L)$, so every subgroup  $H$ of
$G(\Bbb Z/m\Bbb Z)$ with $\pi (H)=L$ is a subgroup of $U$.  Given
$L$ (and hence also $U$) we have the exact sequence $1\to K\to
U\overset{\pi}\to\longrightarrow L\to 1$ and by Lemma 5.1, the
number of $H$ in $U$ with $\pi (H) = L$ is at most $\vert
U\vert^{\ell f(r)}$ where $\ell$ is the derived length of $K$,
$r\leq \dim G$ is the rank of $K$ and $f(r)\leq f(\dim G)$ where
$f$ is some function depending on $r$ and  independent of $m$
(say $f(r)=3r^2+r$).
Now $\vert U\vert \leq m^{\dim G}$ and $K$ being nilpotent, is of
derived length $O(\log\log |K|)$.
 We can, therefore,
deduce that $s_n(G(\Bbb Z/m\Bbb Z))\leq m^{c\dim G f (\dim
G)(\log\log m + \log\dim G)} s_n(G(\Bbb Z/\overline m\Bbb Z))$ for some
constant
$c$  which proves our claim.
\enddemo

Corollary 1.2 shows us that in order to estimate $\alpha_+(G(\Bbb
Z))$ one should concentrate  on $s_n(G(\Bbb Z/m\Bbb Z))$ with
$m\leq n$.  Corollary 5.2 implies that we can further assume that
$m$ is a product of different primes.  So let us now assume that
$m = \prod\limits^t_{i=1} q_i$ where the $q_i$ are different
primes and so $G(\Bbb Z/m\Bbb Z)\simeq \prod G(\Bbb Z/q_i\Bbb Z)$
and $t\leq (1+ o(1))\frac{\log m}{\log\log m}$.  We can further
assume that we are counting only fully proper subgroups of $G(\Bbb
Z/m\Bbb Z)$, i.e., subgroups $H$ which do not contain $G(\Bbb
Z/q_i\Bbb Z)$ for any $1\leq i\leq t$, or equivalently the image
of $H$ under the projection to $G(\Bbb Z/q_i\Bbb Z)$ is a proper
subgroup (see \cite{Lu}). Thus $H$ is contained in $\prod\limits^t_{i=1}
M_i$
where $M_i$ is a maximal subgroup of $G(\Bbb Z/q_i\Bbb Z)$.

Let us now specialize to the case $G=SL_2$, and let $q$ be a
prime.

Maximal subgroups of $SL_2(\Bbb Z/q\Bbb Z)$ are conjugate
 to one of the following three types of subgroups (see \cite{La, Theorem 2.2, 2.3, pp. 183-185}).
\vskip 10pt
{\bf (1)} $B=B_q$ \; -the Borel subgroup of all upper triangular
matrices in $SL_2$.
\vskip 8pt
{\bf (2)} $D=D_q$ \;  -a dihedral subgroup of order $2(q-1)$ or  $2(q+1)$ which is equal
to $N(T_q)$ the normalizer of a split or non-split torus $T_q$.  The group
$T_q$ is either the diagonal subgroup or is obtained as follows: Let $\Bbb F_{q^{2}}$ be the field
of order $q^2, \Bbb F^\times_{q^{2}}$ acts on $\Bbb F_{q^{2}}$ by
multiplication. The latter is a $2$-dimensional vector space over
$\Bbb F_q$.  The elements of norm $1$ in $\Bbb F^\times_{q^{2}}$
induce the subgroup $T_q$ of $SL_2(\Bbb F_q)$.
\vskip 8pt
{\bf (3)} $A=A_q$-a subgroup of $SL_2(\Bbb Z/q\Bbb Z)$ which is of
order at most 120.
\vskip 10pt

There are only boundedly many conjugacy classes of each type.  Also, the number of conjugates of every subgroup
is small, so it suffices to count only subgroups of $SL_2(\Bbb
Z/m\Bbb Z)$ whose projection to $SL_2(\Bbb Z/q\Bbb Z)$ (for
$q\vert m$) is inside either $B, D$, or $A$.

Let $S\subseteq \{q_1\ldots, q_t\}$ be the subset of the prime
divisors of $m$ for which the projection of $H$ is in $A_{q_{i}}$
and $\overline S$ the complement to $S$.  Let $\overline{m} =
\prod\limits_{q\in\overline S} q$ and $\overline H$ the projection
of $H$ to $SL_2(\Bbb Z/\overline{m}\Bbb Z)$.  So $\overline{H}$ is
a subgroup of index at most $n$ in $SL_2(\Bbb Z/\overline{m}\Bbb
Z)$ and the kernel $N$ from $H\to\overline H$ is inside a product
of $\vert S\vert$ groups of type A. As every subgroup of
$SL_2(\Bbb Z/q\Bbb Z)$  is generated by two elements, $H$ is
generated by at most $2 \frac{\log m}{\log\log m}\leq 2\frac{\log
n}{\log\log n}$ generators. Set $k=[2\frac{\log n}{\log\log n} + 1]$ and
choose $k$ generators for $\overline H$.  By a lemma of Gasch\"utz  (cf.
\cite{FJ, Lemma 15.30}) these $k$  generators can be lifted up to
give $k$ generators for $H$.  Each generator can be lifted up in
at most $|N|$ ways and $N$ is a group of order at most
$120^{|S|}\le 120^t\le 120^{\frac{\log n}{\log\log n}}$. We,
therefore, conclude that given $\overline H$ the number of
possibilities for $H$ is at most $120^{2(\log n)^{2}/(\log\log
n)^{2}}$ which is small w.r.t. $n^{\ell (n)}$.

We can, therefore, assume that $S=\phi$ and all the
projections of $H$ are either into groups of type $B$ or $D$.

Now, $B_q$ , the Borel subgroup of $SL_2(\Bbb Z/q\Bbb Z)$,
has a normal unipotent cyclic subgroup $U_q$ of order $q$.  Let
now $S$ be the subset of $\{q_1,\ldots, q_t\}$ for which the
projection is in $B$ and $\overline S$-the complement.  Then
$$H\leq \prod\limits_{q\in S} B_q\times \prod\limits_{q\in\overline
S} D_q.$$  Let $\overline H$ be the projection of $H$ to
$\prod\limits_{q\in S} B_q/U_q\times \prod\limits_{q\in\overline
S} D_q$.  The kernel is a subgroup of the cyclic group
$U=\prod\limits_{q\in S} U_q$.  By Lemma  5.1 we know that given
$\overline H$, there are only few possibilities for $H$. We are,
therefore, led to counting subgroups in 
$$L = \prod\limits_{q\in
S} B_q/U_q\times \prod\limits_{q\in\overline S} D_q.$$  Let $E$
now be the product $$\prod\limits_{q\in S} B_q/U_q\times
\prod\limits_{q\in\overline S} T_q,$$ and for a subgroup $H$ of $L$
we denote $H\cap E$ by $\overline H$.

Our next goal will be to show that given $\overline H$ in $E$,
 the number of possibilities for $H$ is small.  To this end we
formulate first two easy lemmas, which will be used in the proof
of Proposition 5.6 below.  This proposition will complete the
main reduction.

\proclaim{Lemma 5.3}  Let $H$ be a subgroup of $U=U_1\times U_2$.
For $i=1,2$ denote $H_i = \pi_i(H)$ where $\pi_i$ is the
projection from $U$ to $U_i$, and $H^0_i = H\cap U_i$.  Then:
\vskip 10pt
{\bf (i)} $H^0_i$ is normal in $H_i$ and $H_1/H^0_1\simeq H_2/H^0_2$
with an isomorphism $\varphi$ induced by the inclusion of
$H/(H^0_1\times H^0_2)$ as a subdirect product of $H_1/H^0_1$ and
$H_2/H^0_2$,
\vskip 6pt
{\bf (ii)} $H$ is determined by:
\vskip 3pt
\qquad {\bf (a)} $H_i$ for $i=1,2$
\vskip 3pt
\qquad {\bf (b)} $H^0_i$ for $i=1,2$
\vskip 3pt
\qquad {\bf (c)} the isomorphism $\varphi$ from $H_1/H^0_1$
to $H_2/H^0_2$.
\endproclaim

\demo{Proof}  See \cite{Su, p 141}. \qed

\enddemo

\proclaim{Definition 5.4}  Let $U$ be a group and $V$ a subnormal
subgroup of $U$.  We say that $V$ is co-poly-cyclic in $U$ of
co-length $\ell$ if there is a sequence $V=V_0\vartriangleleft
V_1\vartriangleleft\ldots\vartriangleleft V_\ell = U$ such that
$V_i/V_{i-1}$ is cyclic for every $i=1,\ldots ,\ell$.
\endproclaim

\proclaim{Lemma 5.5}  Let $U$ be a group and $F$ a subgroup of
$U$. The number of subnormal co-poly-cyclic subgroups $V$ of $U$
containing $F$ and of co-length $\ell$ is at most $\vert U:F\vert^\ell$.
\endproclaim

\demo{Proof} For $\ell = 1$, $V$ contains $[U,U]F$ and so it
suffices to prove the lemma for the abelian group $\overline U =
U/[U,U]F$ and $\overline F = \{ e\}$.  For an abelian group
$\overline U$, the number of subgroups $V$ with $\overline U/V$
cyclic is equal, by Pontrjagin duality, to the number of cyclic
subgroups.  This is clearly bounded by $\vert\overline U\vert\leq
\vert U:F\vert$.  If $\ell > 1$, then by induction the  number of
possibilities for $V_1$ as in Definition 5.4 is bounded by
$\vert U:F\vert^{\ell -1}$.  Given $V_1$, the number of possibilities for
$V$ is at most $\vert V_1:F\vert\leq \vert U:F\vert$
by the case $\ell = 1$.  Thus,
$V$ has at most $\vert U:F\vert^\ell$ possibilities.  \qed
\enddemo

\proclaim{Proposition 5.6}  Let $D = D_1\times\ldots\times D_s$
where each $D_i$ is a finite dihedral group with a cyclic
subgroup $T_i$ of index $2$.  Let $T=T_1\times\ldots\times T_s$,
so, $\vert D:T\vert = 2^s$.  The number of subgroups $H$ of $D$ whose
intersection with $T$ is a given subgroup $L$ of $T$ is at most
$\vert D\vert^8 2^{2s^{2}}$.
\endproclaim

\demo{Proof}  Denote $F_i = \prod\limits_{j\geq i} D_i$.
We want to count the number of subgroups $H$ of $D$ with
$H\cap T=L$.  Let $L_i =\text{\rm proj}_{F_{i}}(L)$ i.e., the
projection of $L$ to $F_i$, and $\tilde L_{i+1} = L_i\cap
F_{i+1}$, so $\tilde L_{i+1}\subseteq L_{i+1}$.  Let $H_i$ be
the projection of $H$ to $F_i$. Given $H$, the sequence
$(H_1 = H, H_2,\ldots ,H_s)$ is determined and, of course,
vice versa.  We will actually prove that the number of
possibilities for $(H_1,\ldots ,H_s)$ is at most  $\vert
D\vert^8 2^{2s^{2}}$.
\enddemo

Assume now that $H_{i+1}$ is given.  What is the number of
possibilities for $H_i$?  Well, $H_i$ is a subgroup of $F_i =
D_i\times F_{i+1}$ containing $L_i$, whose projection to
$F_{i+1}$ is $H_{i+1}$ and its intersection with $F_{i+1}$,
which we will denote by $X$, contains $\tilde L_{i+1}$.  By
Lemma 4.2, $H_i$ is determined by $H_{i+1} , X,Y,Z$ and
$\varphi$ where $Y$ is the projection of $H_i$ to $D_i$, $Z
= H_i\cap D_i$ and $\varphi$ is an isomorphism from $Y/Z$
to $H_{i+1}/X$.  Now, every subgroup of the dihedral group
is generated by two elements and so the number of
possibilities for $Y$ and $Z$ is at most $\vert D_i\vert^2$
each, and the number of automorphisms of $Y/Z$ is also at
most $\vert D_i\vert^2$.

Let us now look at $X:X$ is a normal subgroup of $H_{i+1}$
with $H_{i+1}/X$ isomorphic to $Y/Z$, so it is meta-cyclic.
Moreover, $X$ contains $\tilde L_{i+1}$.  So by Lemma 4.3,
the number of possibilities for $X$ is at most $\vert H_{i+1}:\tilde
L_{i+1}\vert ^2$.

Now
$\vert H_{i+1}:\tilde L_{i+1}\vert \leq \vert H_{i+1}:L_{i+1}\vert\vert
L_{i+1}:\tilde
L_{i+1}\vert$.  We know that $\vert H_{i+1}:L_{i+1}\vert  =
\vert \text{proj}_{F_{i+1}}(H):\text{proj}_{F_{i+1}}(L)\vert \leq \vert
H:L\vert \leq 2^s$ and
$\vert L_{i+1}:\tilde L_{i+1}\vert =
\vert\text{proj}_{F_{i+1}}(L_i):F_{i+1}\cap L_i\vert\leq \vert
D_i\vert$.  So, $\vert H_{i+1}:\tilde L_{i+1}\vert\leq 2^s\cdot \vert
D_i\vert$.

Altogether, given $H_{i+1}$  (and $L$ and hence also
$L_i$'s and $\tilde L_i$'s) the number of possibilities for
$H_i$ is at most $\vert D_i\vert^8 2^{2s}$.  Arguing, now by
induction we deduce that the number of possibilities for
$(H_1,\ldots ,H_s)$ is at most $\vert D\vert^8 2^{2s^{2}}$
as claimed.  \qed

\vskip 10pt
Let's now get back to $SL_2$:  Proposition 5.6 implies, in the
notations before Lemma 5.3, that when counting
subgroups of $$L =
\prod\limits_{q\in S} B_q/U_q \times\prod\limits_{q\in\overline
S} D_q,$$ we can count instead the subgroups of $$E =
\prod\limits_{q\in S} B_q/U_q\times \prod\limits_{q\in\overline
S} T_q$$ where $T_q$ is a torus in $SL_2(\Bbb Z/q\Bbb Z)$ (so $T_q$ is a
cyclic group of order $q-1$ or $q+1$ while $B_q/U_q$ is a cyclic group of
order $q-1$).

A remark is needed here:  Let $H$ be a subgroup of index at most
$n$ in $SL_2(\Bbb Z/m\Bbb Z)$ which is contained in $X =
\prod\limits_{q\in S} B_q\times \prod\limits_{q\in\overline S}
D_q$ and contains $Y = \prod\limits_{q\in S} U_q\times
v\prod\limits_{q\in\overline S} \{ e\}$.  By our analysis in this
section, these are the groups which we have to count in order to
determine $\alpha_+ (SL_2(\Bbb Z))$.  We proved that for counting
them, it suffices for us to count subgroups of $X_0/Y$ where $X_0
= \prod\limits_{q\in S} B_q\times \prod\limits_{q\in\overline S}
T_q$.  Note though that replacing $H$ with its intersection with
$X_0$, may enlarge the index of $H$ in $SL_2(\Bbb Z/m\Bbb Z)$.
But the factor is at most $$2^{\log m/\log\log m} = m^{1/\log\log
m} \leq n^{1/\log\log n}.$$ As $n\to \infty$, this factor is small
with respect to $n$.  By the remark made in \S1, we can deduce
that our original problem is now completely reduced to the
following extremal problem on counting subgroups of finite
abelian groups: \vskip .20in

Let $\Cal P_- = \{ q_1,\dots, q_t\}$ and $\Cal P_+=\{ q'_1,\dots,
q'_{t'}\}$ be two sets of (different) primes and let $\Cal P = \Cal
P_-\bigcup \Cal P_+$.  Denote
$$
f(n)=\sup\{s_r(X)|X=\prod\limits^t_{i=1} C_{q_i - 1} \times
\prod\limits^{t'}_{i=1} C_{q'_i +1}\}$$ where the supremum is
over all possible choices of $\Cal P_-, \Cal P_+$ and $r$ such
that $$r\prod\limits^t_{i=1}q_i\prod\limits^{t'}_{j=1}q'_j \; \le \; n,$$
and where $C_m$ denotes the cyclic group of order $m$. The discussion above implies:

\proclaim{Proposition 5.7} We have
$$\alpha_+ (SL_2(\Bbb Z)) = \overline{\lim} \;\frac{\log f(n)}{\lambda (n)}.
$$
\endproclaim

\vskip .1in
{\bf \S6. Counting subgroups of $p$-groups}
\vskip .10in

In this section we first give some general estimates for the
number of subgroups of finite abelian $p$-groups which will be
needed in \S 7. As an application we obtain a lower bound
for the subgroup growth of uniform pro-$p$-groups (see definitions
later).

For an abelian $p$-group $G$, we denote by $\Omega_i(G)$
 the subgroup of elements of order dividing $p^i$.  Then
$\Omega_i(G)/\Omega_{i-1}(G)$ is an elementary abelian
group of order say $p^{\lambda_{i}}$ called the $i$-th {\it
layer} of $G$.  We call the sequence
$\lambda_1\geq\lambda_2\geq\ldots\geq \lambda_r$ the
{\it layer type} of $G$. It is clear that this sequence is
decreasing.

Denote by $\bmatrix \lambda \\ \nu\endbmatrix_p$ the
$p$-binomial coefficient, that is, the number of
$\nu$-dimensional subspaces of a $\lambda$-dimensional vector
space over $\Bbb Z/p\Bbb Z$.

The following holds (see \cite{LS, Proposition 1.5.2}).

\proclaim{Proposition 6.1}
\vskip 8pt
{\bf (i)} $p^{\nu (\lambda -\nu )} \leq \bmatrix \lambda \\
\nu\endbmatrix_p \leq p^\nu\cdot p^{\nu (\lambda -\nu
)}$.
\vskip 3pt
{\bf (ii)} $\max\bmatrix \lambda \\ \nu\endbmatrix_p$ is
attained for $\nu = [\frac{\lambda}{2}]$ in which case
$\bmatrix \lambda \\ \nu\endbmatrix_p = p^{\frac{1}{4}
\lambda^{2} + O(\lambda)}$ holds as $\lambda\to\infty$.
\endproclaim

The starting point is the following well-known formula (see\cite{Bu}).

\proclaim{Proposition 6.2}  Let $G$ be an abelian $p$-group of
layer type $\lambda_1\geq\lambda_2\ldots \geq\lambda_r$. The number of
subgroups  of layer type $\nu_1\geq \nu_2\ldots$ is
$$
\prod\limits_{i\geq 1} p^{\nu_{i+1}(\lambda_{i}-\nu_{i})}
\bmatrix \lambda_i - \nu_{i+1} \\ \nu_i -
\nu_{i+1}\endbmatrix_p.  \qquad\qquad\qquad\qed
$$
\endproclaim

(In the above expression we allow some of the $\nu_i$ to be $0$.)

We need the following estimate.

\proclaim{Proposition 6.3}
$$
\prod\limits_{i\geq 1} p^{\nu_{i}(\lambda_{i}-\nu_{i})}
\leq \prod\limits_{i\geq 1}
p^{\nu_{i+1}(\lambda_{i}-\nu_{i})} \bmatrix
\lambda_i-\nu_{i+1} \\ \nu_{i} - \nu_{i+1}\endbmatrix_p
\leq p^{\nu_{1}} \prod\limits_{i\geq 1}
p^{\nu_{i}(\lambda_{i}-\nu_{i})} .
$$
\endproclaim

\demo{Proof}  By Proposition 6.1 we have
$$\align
\prod\limits_{i\geq 1} p^{\nu_{i+1}(\lambda_{i}-\nu_{i})}
\bmatrix \lambda_i-\nu_{i+1}  \\ \nu_{i} - \nu_{i+1}\endbmatrix_p
\; &\leq  \;\, \prod\limits_{i\geq 1}
p^{\nu_{i+1}(\lambda_{i}-\nu_{i})}
\cdot
p^{(\nu_{i}-\nu_{i+1})((\lambda_{i}-\nu_{i+1})-
(\nu_i-\nu_{i+1}))}
\cdot
p^{(\nu_{i}-\nu_{i+1})} \\
&= \;\,  p^{\nu_{1}}
\prod\limits_{i\geq 1}
p^{\nu_{i+1}(\lambda_{i}-\nu_{i})}\cdot
p^{(\nu_{i}-\nu_{i+1})(\lambda_{i}-\nu_{i})} = p^{\nu_{1}}
\prod\limits_{i\geq 1} p^{\nu_{i}(\lambda_{i}-\nu_{i})}.
\endalign
$$
The lower bound follows in a similar way. \qed
\enddemo

\proclaim{Corollary 6.4}   Let $G$ be an abelian group of order
$p^\alpha$ and layer type
$\lambda_1\geq\lambda_2\geq\ldots\geq\lambda_r$. Then $\vert
G\vert^{-1}\prod\limits_{i\geq 1} p^{\lambda_{i}^{2}/4} \; \leq \;
\vert\text{\rm Sub}(G)\vert  \; \leq  \; \vert G\vert^2 \prod\limits_{i\geq 1}
p^{\lambda_{i}^{2}/4}$ holds.
\endproclaim

\demo{Proof}  Considering subgroups $H$ of layer type
$[\frac{\lambda_{1}}{2}]\geq
[\frac{\lambda_{2}}{2}]\geq\ldots$ we obtain that \newline
$\vert\text{\rm Sub}(G)\vert\geq\prod\limits_{i\geq 1}
p^{[\frac{\lambda_{i}}{2}](\lambda_{i}-[\frac{\lambda_{i}}{2}])}
\geq p^{-r} \prod\limits_{i\geq 1}
p^{\lambda_{i}^{2}/4}$ which implies the lower bound.
\enddemo

On the other hand, for any fixed layer type
$\nu_1\geq\nu_2\geq\ldots $ the number of subgroups $H$
with this layer type is at most
$$
p^{\nu_{1}} \prod\limits_{i\geq 1} p^{\nu_{i}
(\lambda_{i}-\nu_{i})}  \; \leq  \; \vert G\vert \prod\limits_{i\geq
1} p^{\lambda_{i}^{2}/4}.
$$

The number of possible layer types
$\nu_1\geq\nu_2\geq\ldots$ of subgroups of $G$ is
bounded by the number of partitions of the number
$\alpha$ hence it is at most $2^\alpha\leq \vert G\vert$.
This implies our statement. \qed

Let us make an amusing remark which will not be needed later.

If $G$ is an abelian $p$-group of the form
$G = C_{x_{1}}\times C_{x_{2}}\times\ldots\times
C_{x_{t}}$ then it is known (see \cite{LS,
\S 1.10}) that
$\vert\text{End} (G)\vert = \prod\limits_{j,k\geq 1} \text{ gcd} (x_{j},
x_k)$.  Noting that $\prod\limits_{j,k\geq 1} \text{gcd}(x_j,x_k) =
\prod\limits_{i\geq 1} p^{\lambda^{2}_{i}}$ we obtain that
$$
\vert G\vert^{-1} \vert\text{End}(G)\vert^{\frac{1}{4}} \; \leq \;
\vert\text{Sub}(G)\vert   \;  \leq   \;  \vert G\vert^2
\vert\text{End}(G)\vert^{\frac{1}{4}}.
$$
These inequalities clearly extend to arbitrary finite abelian
groups $G$.

For the application of the above results to estimating the subgroup growth
of $SL_d(\Bbb Z _p)$ we have to introduce additional notation. For a group
$G$ let $G^k$ denote the subgroup generated by all $k$-th powers. For odd
$p$ a {\it powerful} $p$-group $G$ is a $p$-group with the property that
$G/G^p$ is abelian. (In the rest of this section we will always assume that
$p$ is odd,the case $p=2$ requires only slight modifications.) $G$ is
said to be {\it uniformly powerful} ({\it uniform}, for short) if it is
powerful and the indices $\vert G^{p^i}:G^{p^{i+1}}\vert$ do not depend on
$i$ as long as $i<e$, where $p^e$ is the exponent of $G$.

Now let $G$ be a uniform group of exponent $p^e$ , where $e=2i$ , with $d$
generators. Then $G^{p^i}$ is a homocyclic abelian group of exponent $p^i$
and $d$ generators (i.e. it has layer type $d,d,\ldots,d$ with $i$ terms)
\cite{Sh}.

Consider subgroups $H$ of $G^{p^i}$ of layer type
$\nu,\nu,\ldots,\nu$ ($i$ terms). The number of such subgroups is at least
$p^{i\nu(d-\nu)}$ by Proposition 6.3. . The index $n$ of such a subgroup
$H$
in $G$ is $p^{di+(d-\nu)i}$. Hence the number of index $n$ subgroups in
$G$
is at least $n^x$ where $x=\frac{\nu(d-\nu)}{2d-\nu}$. Substituting
$\nu= [d(2-\sqrt{2})]$ we see that $x$ can be as large as
$(3-2\sqrt{2})d-(\sqrt{2}-1)$.

Let now $U$ be a uniform pro-$p$-group of rank $d$ , i.e. an inverse limit
of $d$-generated finite uniform groups $G$. Then we see that for infinitely
many $n$ we have $s_n(G)\;\geq \;n^{(3-2\sqrt{2})d-(\sqrt{2}-1)}$.

Now $SL_d(\Bbb Z _p)$ is known to have a finite index uniform pro-$p$-
subgroup of rank $d^2-1$ 
 (see\cite{DDMS, Theorem 5.2}). This proves the following

\proclaim{Proposition 6.5} The group $SL_d(\Bbb Z _p)$ has subgroup growth of type at
least   $n^{(3-2\sqrt{2})d^2-2(2-\sqrt{2})}.$
\endproclaim

B. Klopsch proved \cite{Kl} that if $G$ is a residually finite virtually
soluble minimax group of Hirsch length $h(G)$ then its subgroup growth is
of type at least $n^{h(G)/7}$. By using the above argument one can improve
this to $n^{(3-2\sqrt{2})h(G)-(\sqrt{2}-1)}$.

 \vskip .2in
{\bf \S7. Counting subgroups of abelian groups}
\vskip .10in

The aim of this section is to solve a somewhat unusual
extremal problem concerning the number of subgroups of
abelian groups.  The result we prove is the crucial ingredient
in obtaining a sharp upper bound for the number of
congruence subgroups of $SL(2,\Bbb Z)$. Actually we prove a slightly more general result which will be used in [LN] to obtain similar bounds for  other arithmetic groups.

We will use Propositions 6.2 and 6.3 in conjunction with the
following simple (but somewhat technical) observations.

\pagebreak

\proclaim {Proposition 7.1} Let $R\geq 1$ and let $C,t\in \Bbb N$ be fixed. Consider pairs of sequences $\{\lambda_i\},\{\nu_i\}$ of nonnegative integers, such that $\lambda_i\leq t$ for all $i$ and $\sum_{i\geq 1}(R\lambda_i+\nu_i) \leq C$.

Under these conditions the maximal value of the expression
$A(\{\lambda\},\{\nu\})=\sum_{i\geq 1}\nu_i(\lambda_i-\nu_i)$ can be attained by a pair of sequences $\{\lambda_i\},\{\nu_i\}, \ i=1,2,..,r$ such that: \medskip

{\bf (i)} $\lambda_1\geq \lambda_2 \geq ... \geq \lambda_r, \quad \nu_1\geq \nu_2 \geq ... \geq \nu_r\geq 1$, and $\lambda_i \geq \nu_i$ for all $i$, \medskip

{\bf (ii)} $\lambda_1=\lambda_2 = ... = \lambda_{r-1}=t$ and \medskip

{\bf (iii)} for some $0\leq b\leq r-1$ we have
$\nu_1= \nu_2 = ... = \nu_b=\nu_{b+1}+1=...=\nu_{r-1}+1.$
If $\lambda_r=t$ then also $\nu_r \in \{\nu_1, \nu_1-1\}$.
\endproclaim

\demo{Proof} Suppose the maximum of $A(\{\lambda\},\{\nu\})$ is attained by a pair $\{\lambda_i\},\{\nu_i\}$ of sequences of non-negative integers. Deleting pairs with $\nu_j=0$ does not change the value of $A(\{\lambda\},\{\nu\})$
hence we can assume that all $\nu_i \geq 1$. If $\lambda_j < \nu_j$ for some $j$, then we can delete $\lambda_j$ and $\nu_j$ from the sequences and in this way the value of $A(\{\lambda\},\{\nu\})$ increases, a contradiction. Hence we have that $\lambda_i \geq \nu_i$ for all $i$. By relabelling the indices we can further assume that $\nu_1\geq \nu_2 \geq ... \geq \nu_r \geq 1$. \medskip

Now, if $\pi$ is a permutation of $\{1,2,...,r\}$, it is clear that the maximum of $\sum_i\lambda_{\pi(i)}\nu_i$ (and hence of $A(\{\lambda_{\pi(i)}\},\{\nu_i\})$) is achieved for permutations $\pi$ such that $\lambda_{\pi(1)}\geq \lambda_{\pi(2)} \geq ... \geq \lambda_{\pi(r)}$. By the maximality of the pair $\{\lambda_i\},\{\nu_i\}$ it now follows that
$\lambda_1\geq \lambda_2 \geq ... \geq \lambda_r$ as well, proving (i). We shall call a pair of sequences $\{\lambda\},\{\nu\}$ satisfying (i) {\it good}. \medskip 

Let $j$ be the smallest index
such that we have $t > \lambda_j\geq\lambda_{j+1} \geq 1$
(if there is no such $j$ then (ii) holds).

Assume that $\lambda_{j+1}=\ldots =\lambda_{j+k}$ and
$\lambda_{j+k}>\lambda_{j+k+1}$ or $j+k=r$. The
condition $\nu_j\geq\nu_{j+k}$ implies that
$\nu_j((\lambda_j + 1)-\nu_j) + \nu_{j+k}
((\lambda_{j+k}-1)-\nu_{j+k}) \geq \nu_j
(\lambda_j-\nu_j) + \nu_{j+k} (\lambda_{j+k} - \nu_{j+k}).
$
If $\lambda_{j+k} = \nu_{j+k}$ then (by deleting some terms and
relabelling the rest) we can replace
our sequences by another good pair for which $\sum\limits_{i\geq
1} \lambda_j$ is strictly smaller and the value of $A(\{\lambda_i\}, \{\nu_i\})$ is the same.  Otherwise, replacing
$\lambda_j$ by $\lambda_j +1$ and $\lambda_{j+k}$ by
$\lambda_{j+k} -1$ we obtain a good pair of sequences for
which $\{\lambda_i\}$ is lexicographically strictly greater
and for which $A(\{\lambda_i\}, \{\nu_i\})$ is at least as large (hence maximal).
\enddemo

It is clear that by repeating these two types of moves we
eventually obtain a good pair $\{\lambda_i\},\{\nu_i\}$
satisfying (ii) as well.
\medskip

Now set $\beta = \nu_1 + \nu_2 +\ldots +\nu_{r-1}$.  Then
$$
\sum\limits_{i\geq 1} \nu_i(\lambda_i-\nu_i) = t\beta -
(\nu^2_1 +\ldots +\nu^2_{r-1}) + \nu_{r} (\lambda_{r} -
\nu_{r}).
$$

It is clear that if the value of such an expression is maximal,
then the difference of any two of the $\nu_j$ with $j\leq r-1$
is at most $1$.  Part (iii) follows. \qed
\bigskip

\proclaim{Proposition 7.2}  Let $x_1,x_2,\dots,x_t$ be positive
integers such that at most $d$ of the $x_i$ can be equal. Then
$$
\prod_{i=1}^t x_i \; \geq  \; \bigg(\frac{t}{ed}\bigg)^t$$
holds.
\endproclaim

\demo{Proof} If say, $x_1$ is the largest among the $x_i$ then
$x_1\ge\frac{t}{d}$. By induction we can assume that
$\prod\limits_{i=2}^t x_i\geq\bigg(\frac{t-1}{ed}\bigg)^{t-1}$ holds. Then
$$
\prod_{i=1}^tx_i\geq\frac{t}{d}\bigg(\frac{t-1}{ed}\bigg)^{t-1}\geq
e\bigg(\frac{t}{ed}\bigg)\bigg(\frac{t-1}{ed}\bigg)^{t-1}\geq
e  \bigg(\frac{t}{ed}\bigg)^t\bigg(\frac{t-1}{t}\bigg)^{t-1}=
$$
$$
=\bigg(    \frac{t}{ed}\bigg)^t
\frac{e}{\bigg(1+\frac{1}{t-1}\bigg)^{t-1}}\geq
\bigg(\frac{t}{ed}\bigg)^t,\ \ \ \ \text{ }\hskip2cm \qed
$$
as required.
\enddemo

The main result of this section is the following.

\proclaim{Theorem 7.3}  Let $R\geq 1$ be a real number and $d$ be a fixed integer $\geq 1$.  Let
$n,r$ be positive integers.  Let $G$ be an abelian group of the
form $G = C_{x_{1}}\times C_{x_{2}}\times\ldots\times C_{x_{t}}$
where at most $d$ of the $x_i$ can be  equal. Suppose that
$r\vert G\vert^R\leq n$ holds.  Then the number of subgroups  of
order $\leq r$ in $G$ is at most $n^{(\gamma + o(1))\ell (n)}$
where $\gamma = \frac{(\sqrt{R(R+1)}-R)^2}{4R^2}$. In particular if $R=1$ then $\gamma =\frac{3-2\sqrt{2}}{4}$.
\endproclaim

\demo{Proof}  We start the proof with several claims.
\enddemo

{\bf Claim 1.}  $t\leq (1 + o(1)) \ell (n)$.

\demo{Proof}  By Proposition 7.2 we have
$\big(\frac{t}{ed}\big)^t\leq n$. This easily implies the claim.
\enddemo

{\bf Claim 2.}  In proving the theorem, we may assume that
$t\geq \gamma \ell (n)$.

\demo{Proof}  For otherwise, every subgroup of $G$ can be
generated by $\gamma \ell (n)$ elements hence
$\vert\text{\rm Sub}(G)\vert\leq \vert G\vert^{\gamma \ell (n)}
\leq n^{\gamma \ell (n)}$.
\enddemo

Now let $a(n)$ be a monotone increasing function which
goes to infinity sufficiently slowly.  For example, we may set
$a(n) = \log\log\log\log n$.

Let $G_p$ denote the Sylow $p$-subgroup of $G$ and let
$\lambda_1^p \geq \lambda_2^p\geq\ldots$ denote the
layer type of $G_p$. Loosely speaking, we call any layer of some $G_p,$ a layer of G. 
 We call such a layer {\it essential} if
its dimension $\lambda_i^p$ is at least
$\frac{\ell(n)}{a(n)}$.  Clearly the essential layers in $G_p$ correspond
to the layers of a certain subgroup $E_p$ of $G_p$ (which
equals $\Omega_i(G_p)$ for the largest $i$ such that
$\lambda_i^p\geq \frac{\ell(n)}{a(n)}$).  Let us call $E =
\prod\limits_p E_p$ the {\it essential subgroup} of $G$.

{\bf Claim 3.}  Given $E\cap T$ we have at most $n^{o(\ell
(n)}$ (i.e., a small number) of choices for a subgroup $T$ of $G$.

\demo{Proof}  It is clear from the definitions that every
subgroup of the quotient groups $G_p/E_p$ and hence of
$G/E$ can be generated by
less than $\frac{\ell (n)}{a(n)}$ elements.  Therefore the same is
true for $T/T\cap E$.  This implies the claim.
\enddemo

By Claim 3, in proving the theorem, it is sufficient to
consider subgroups $T$ of $E$.

Let $v$ denote the exponent of $E$. Then $E$ is the subgroup
of elements of order dividing $v$ in $G$.
Now $v$ is the product of the exponents of the $E_p$ hence the product
of the exponents of the essential layers of $G$.
It is clear from the
definitions that we have $v^{\ell (n)/a(n)}\leq n$,
hence $v\leq (\log n)^{a(n)}$.  Using well-known estimates of
number theory \cite{Ra}
we immediately obtain the
following.
\vskip 3pt
{\bf Claim 4.}  {\bf  (i)} the number $z$ of different primes
dividing $v$ is at most
$\frac{\log v}{\log\log v}\leq\frac{a(n)\log\log n}{\log\log\log n}$.
\vskip 3pt
{\bf (ii)} The total number of divisors of $v$ is at most
$v^{\frac{c}{\log\log v}}\leq\log n^{\frac{ca(n)}{\log\log\log n}}$
for some constant $c > 0$.
\vskip 3pt
{\bf Claim 5.}  $\vert G:E\vert\geq (\log n)^{(1+o(1))t}$.

\demo{Proof}  Consider the subgroup
$E^i = E\cap C_{x_i}$.  It follows that $E^i$ is the subgroup
of elements of order dividing $v$ in $C_{x_{i}}$.  Set  $e_i =
\vert E^i\vert$ and $h_i = x_i/e_i$.  It is easy to see that $E
= \prod\limits_{i\geq 1} E^i$, hence $\vert G:E\vert =
\prod\limits_{i\geq 1} h_i$.
\enddemo

By Claim 4(ii) for the number $s$ of different values of the
numbers $e_i$ we have $s = (\log n)^{o(1)}$.  We put the
numbers $x_i$ into $s$ blocks
according to the value of $e_i$.
By our condition on the $x_i$ it follows that at most $d$ of
the numbers $h_i$ corresponding to a given block are equal.
Hence altogether $ds$ of the $h_i$ can be equal. Using
Proposition 7.2 we obtain that
$\vert G:E\vert\geq\prod\limits_{i\geq 1}h_i\geq
\big(\frac{t}{eds}\big)^t$.

Since $sd = (\log n)^{o(1)}$ and by Claim 2 \
$t\geq \gamma \frac{\log n}{\log\log n}$ we obtain that
$\vert G:E\vert\geq (\log n)^{(1+o(1))t}$ as required.

Let us now choose a group $G$ and a number $r$ as in the
theorem for which the number of subgroups $T\leq E$ of
order dividing $r$ is maximal.  To complete the proof it is
clearly sufficient to show that this number is at most
$n^{(\gamma +o(1))\ell (n)}$.

Denote the order of the corresponding essential subgroup
$E$ by $f$ and the index $\vert G:E\vert$ by $m$.

Using Propositions 6.2 and 6.3 we see that apart from an
$n^{o(\ell (n))}$ factor (which we ignore) the number of
subgroups $T$ as above is at most
$$
\prod\limits_{p|f}\prod\limits_{i\geq 1}
p^{\nu^{p}_{i}(\lambda^{p}_{i}-\nu^{p}_{i})}\tag 7.1
$$
for some
$\nu^p_i, \lambda^p_i$ where $\{\lambda_i^p\},\{\nu_i^p\}$ is a
pair of sequences for every $p$,
$\prod\limits_{p}\prod\limits_{i\geq 1}
p^{\lambda^{p}_{i}}$ divides $f$ and
$\prod\limits_{p}\prod\limits_{i\geq 1} p^{\nu^{p}_{i}}$
divides $r$.  Assuming that $f^Rr$
is fixed together with the upper bound $t$ for all the
$\lambda_i^p,\mu_i^p$, let us estimate the value of the expression
(7.1).

By Proposition 7.1 a maximal value of an expression like (7.1) is attained
for a choice of the $\lambda_i^p,\nu_i^p$ (for the sake of simplicity
we use the same notation for the new sequences) such that
for every $p$ there are at most 3 different pairs
$(p^{\lambda_i^p},p^{\nu_i^p})$ equal to say
$$
(p^t,p^{\mu^p+1}), \;\;\, (p^t,p^{\mu^p}),\ \ \text{ and }\ \
   (p^{\tau^p}, p^{\mu_0^p})
$$
where $\mu_0^p \leq \tau^p<t$ and
$\mu^p <t$ for all $p$.

Exchange the pairs equal to the first type for pairs
equal to $(p^t,p^{\mu^p})$. We obtain an expression
like (7.1)  such that the ratio of the two expressions
is at most
$$
\prod_p\prod_{i\geq 1}p^{\lambda_i^p}\leq n\ .
$$

If now there are say $\alpha\,^p$ pairs with $(p^{\lambda_i^p},p^{\nu_i^p})$
equal to $(p^t,p^{\mu^p})$ then take $\beta^p$ to be the largest integer
with $2^{\beta^p}\leq p^{\alpha^p}$ and set $\beta_1^p=
\big[\log_2p\big]$. (Note that for every $p$ there is at most one
pair of the form $(p^{\tau^p}, p^{\mu_0^p})$.)

Consider the expression
$$
\prod_{p}2^{\beta^p\mu^p(t-\mu^p)}\
    2^{\beta_1^p\mu_0^p(\tau^p-\mu_0^p)}.\tag 7.2
$$

Its value may be less than that of (7.1) but in this case their ratio is bounded
by $(2^{2z})^{t^2}n$ (where $z$ is the number of primes
dividing $v$). Hence this ratio is at most
$$
2^{(2+o(1))\ell(n)^2\
  \frac{a(n)\log\log n}{\log\log\log n}}\leq
n^{(2+o(1))\ell(n)\
  \frac{a(n)}{\log\log\log n}}=n^{o(\ell(n))}.
$$
To prove our theorem it is sufficient to bound the value
of (7.2) by $n^{(\gamma+o(1))\ell(n)}$.

It is clear that the value of (7.2) is equal to the value of another
expression
$$
\prod_{k\geq 1}2^{\nu_k(\lambda_k-\nu_k)}    \tag 7.3
$$ for appropriate sequences $\{\lambda_k\}, \{\mu_k\}$
which both have $\sum\limits_p(\beta^p+\beta^p_1)$ terms and for which $\lambda_k,\mu_k \leq t$ and also
$\prod\limits_{k\geq 1}2^{R\lambda_k+\nu_k}\leq f^R\cdot r$, 
i.e.
$\sum_{k \geq 1}(R\lambda_k+\nu_k) \leq \log (rf^R). \hskip 2cm (*)$

More precisely, the sequence $\{\lambda_k\}$ has $\sum\limits_p \beta^p$ terms equal to $t$ and $\beta^p_1$ terms equal to $\tau^p$ for every $p$, while $\{\mu_k\}$ consists of $\mu^p$ repeated $\beta^p$ times and $\mu^p_0$ repeated $\beta^p_1$ times each (in the appropriate order).

By Proposition 7.1 the expression 7.3 attains its maximal
value for some sequences $\{\lambda _k\}, \{\nu _k\}$
such that all but one of the $\lambda _k$, say $\lambda _{a+1}$
are equal to $t$ and we have \newline $\nu_1 =
\nu_2=\ldots = \nu_b = 1 + \nu_{b+1}=\ldots=
1+\nu_a$ for some $b\leq a$.

Consider now the expression
$$
\prod\limits_{k\geq 1} 2^{\nu'_k(\lambda'_{k}-\nu'_{k})}
\tag 7.4
$$
where
$$
t = \lambda'_1 = \ldots = \lambda'_a \quad
(\lambda'_{a+1} = 0)
$$
and $\nu_a = \nu'_1 = \nu'_2 =\ldots =\nu'_a\quad (\nu'_{a+1}=0)$.

It easily follows that the value of (7.3) is at most
$2^{2t^{2}}$ times as large as the value of (7.4) and
$2^{2t^{2}} = n^{o(\ell (n))}$.  Hence it suffices to bound the
value of (7.4) by $n^{(\gamma +o(n))\ell (n)}$.

To obtain our final estimate denote $2^a$ by $y$, $ m^{1/t}$ by
$w$ (where $m=\vert G:E\vert$) and set $x=y\cdot w$.

For some constants between $0$ and $1$ we have
$y = x^\rho$ and $\nu'_1 = \sigma t$.  Then \newline $w =
x^{1-\rho} = y^{\frac{1-\rho}{\rho}}$.

Note that the condition (*) implies $2^{at(R+\sigma)} = y^{\sigma t}y^{Rt} \leq rf^R$.
We have $n\geq r(mf)^R \geq y^{\sigma t} \cdot y^{Rt} \cdot w^{Rt}$
hence $\log n\geq t\cdot\log y \left(R+\sigma
+R\frac{1-\rho}{\rho}\right)$.

By Claim 5 we have $w\geq (\log n)^{(1+ o(1))}$. Hence
$$(1+ o(1))\log\log n \leq \log w = \frac{1-\rho}{\rho} \log
y.$$  Therefore
$$\align \frac{(\log n)^{2}}{\log\log n} &\geq \frac{t^{2} (\log y)^{2}
(R+\sigma +R\frac{1-\rho}{\rho})^{2}}{(\frac{1-\rho}{\rho}
\log y)}\cdot \big(1+ o(1)\big)\\
&= \big(1+ o(1)\big)\cdot t^2 \log y \left(R + \sigma +
R\frac{1-\rho}{\rho}\right)^{2}\cdot \left(\frac{\rho}{1-\rho}\right).\endalign
$$
The value of (7.4) is $y^{\sigma t(t-\sigma t)}$ which as we
saw is an upper bound for the number of subgroups $R$
(ignoring an $n^{o(\ell (n))}$ factor).  Hence

$$\align
  & \frac{\log \text{(number of subgroups $T$)}}
{(\frac{(\log n)^{2}}{\log\log n})} \\  & 
\leq (1+ o(1))\frac{t^{2} \sigma (1-\sigma )\log y}{t^{2}\log y (R+\sigma
+R\frac{1-\rho}{\rho})^{2}(\frac{\rho}{1-\rho})}\\ &=(1+
o(1)) \frac{\sigma (1-\sigma
)(\frac{1-\rho}{\rho})}{(R+\sigma +R\frac{1-\rho}{\rho})^{2}}
= (1+ o(1))\frac{\sigma (1-\sigma
)\rho(1-\rho)}{(R+\rho\sigma )^{2}}.
\endalign$$
As observed in \S 3, the maximum value of
$\frac{\sigma (1-\sigma )\rho(1-\rho)}{(R+\rho\sigma
)^{2}}$ for $\sigma, \rho \in (0,1)$ is $\gamma$.  The proof of the theorem is complete.
\qed

By using a similar but simpler argument, one can also show
the following

\proclaim{Proposition 7.4}  Let $G$ be an abelian group of order
$n$ of the form \newline $G = C_{x_{1}}\times
C_{x_{2}}\times\ldots\times C_{x_{t}}$ where $x_1 > x_2 > \ldots
x_t$.  Then $\vert\text{\rm Sub} (G)\vert \leq n^{(\frac{1}{16 } +
o(1))\ell (n)}$. This bound is attained if $x_{i} = t\cdot i$ for
all $i$.
\endproclaim

Combining this result with an earlier remark, we obtain that
$n^{(\frac{1}{4} + o(1))\ell (n)}$ is the maximal value of
$\prod\limits_{i,j} \text{gcd} (x_i, x_j)$ where the $x_i$ are
different numbers whose product is at most $n$.

Note that $\vert\text{\rm Sub}(G)\vert$ is essentially the number of
subgroups $T$ of order $[\sqrt{\vert G\vert}]$ (see \cite{Bu}
for a strong version of this assertion).  Hence Proposition 7.4
corresponds to the case\break  $R = 1,\; r\sim n^{1/3}$ of Theorem 7.3.

\vskip .20in
{\bf \S8. End of proofs of Theorems 2, 3, and 4.}
\vskip .1in

Theorem 2 is actually proved now:  the lower bound was shown as a
special case of $R=R(G) = 1$ in \S 3.  For  the  upper bound, we have
shown in Proposition 5.7 how $\a_+(SL_2(\Bbb Z))$ is equal to
$\overline{\lim} \frac{\log f(n)}{\lambda(n)}$ (see Proposition 5.7 for the
definition of $f(n)$). But Theorem 7.3 implies, in particular,
that $f(n)$ is at most $n^{(\gamma+o(1)) \ell(n)}$ where $\gamma =
\frac{3-2\sqrt{2}}{4}$.  This proves that $\a_+ (SL_2(\Bbb Z))\le
\gamma$ and finishes the proof.

The proof of Theorem 3 is similar, but several remarks should be
made:  The lower bound was deduced in \S 4.  For the upper bound,
one should follow the reductions made in \S 6.  The proof can be
carried out in a similar way for $SL_2(\Cal O)$ instead of
$SL_2(\Bbb Z)$ but the following points require careful
consideration.

{\bf 1)}  One can pass to the case that $m$ is an ideal which is a
product of different primes $\pi_i$'s in ${\Cal  O}$, but it is
possible that $\Cal  O/\pi_i$ is isomorphic to $\Cal O/\pi_j$.
Still, each such isomorphism   class of quotient fields can occur
at most $d$ times when $d=[k:\Bbb Q]$.

{\bf 2)}  The maximal subgroups of $SL_2(\Bbb F_q)$ when $\Bbb F_q$ is a
finite field of order $q$ ($q$ is a prime power,  not necessarily
a prime) are the same $B, D$ and $A$ as described in (1), (2), and
(3) of \S 5.

The rest of the reduction can be carried out in a similar way to
\S 5.  The final outcome is not exactly as $f(n)$ at the end of \S
5, but can be reduced to a similar problem when $\tilde f(n)$
counts $s_r(X)$ when $X$ is a product of abelian cyclic groups,
with a bounded multiplicity.  Theorem 7.3 covers also this case
and gives a bound to $\tilde f(n)$ which is the same as for
$f(n)$.  Thus $\a_+ (SL_2(\Cal  O)) \le \gamma =
\frac{3-2\sqrt{2}}{4}$.

We finally mention the easy fact, that replacing $\Cal  O$ by
$\Cal  O_S$ when $S$ is a finite set of primes (see the
introduction) does not change $\a_+$ or $\a_-$. To see this one
can  use the fact that for every completion at a simple
prime $\pi$ of $\Cal  O$, $G(\Cal  O_\pi)$ has polynomial subgroup
growth and then use the well known techniques of subgroup growth
and the fact that
$$
G(\hat \Cal  O) = G(\hat{\Cal  O}_S)
\times\mathop{\pi}\limits_{\pi \in S\setminus V_\infty} G({\Cal
O}_\pi)
$$
to deduce that $\a (G(\hat{\Cal  O})) = \a (G(\hat {\Cal  O}_S))$.

Another way to see it, is to observe that $G(\hat{\Cal  O}_S)$
is a quotient of $G(\hat{\Cal  O})$, and, hence,  $\a_+ (G({\Cal  O}))
\ge \a_+(G({\Cal  O}_S))$.  On the other hand, the proof of the
lower bound for $\a (G({\Cal  O}))$ clearly works for $G({\Cal
O}_S)$. Theorem 3 is, therefore, now proved, as well as Theorem 4 (since we have not
used the GRH for the upper bounds in Theorem 3).

\vskip .2in
{\bf \S9. An extremal problem in elementary number theory.}
\vskip .10in

  The counting techniques in this paper can be applied to solve the following
extremal problem in multiplicative number theory.

For $n \to \infty$, let
$$\align  &M_1(n) = \max \Bigg\{ \prod_{1 \le i, j \le t} \text{gcd}(a_i, a_j) \;\,
\Bigg | \;\, 0 < t, a_1 <  a_2 <  \ldots <  a_t \in \Bbb Z,
\; \, \prod_{i=1}^t a_i \le n\Bigg\},\\
&M_2(n) = \max \Bigg\{ \prod_{p, p' \in \Cal P} \text{gcd}(p-1,\;
p'-1) \;\, \Bigg | \;\, \Cal P = \text{set of distinct primes
where} \, \prod_{p\in\Cal P} p \le n\Bigg\}.\endalign$$ We shall
prove the following theorem which can be considered as a baby
version of Theorem 2 (compare also to Theorem 7.3 ). Note that Theorem 9.1 immediately
implies
Theorem 9. \proclaim{Theorem 9.1}   Let $\lambda(n) = \frac{(\log
n)^2}{\log\log n}.$ Then
$$\underline{\lim} \; \frac{\log M_1(n)}{\lambda(n)} = \overline{\lim} \; \frac{\log M_2(n)}{\lambda(n)} = \frac14.$$
\endproclaim

\demo{Proof} Recall that if $a_1, a_2, \ldots, a_t \in \Bbb Z$ and $G = C_{a_1}\times C_{a_2} \times \cdots \times
C_{a_t}$ is a direct product of cyclic groups then by \S7,
$$|G|^{-1} |\text{End}(G)|^{\frac14} \; \le \; |\text{Sub}(G)| \; \le \; |G|^2  \; |\text{End}(G)|^{\frac14},$$
and
$$|\text{End}(G)| \; = \; \prod_{1\le i, j\le t} \text{gcd}(a_i, a_j).$$
Proposition 7.4 implies that
$$\overline{\lim} \; \frac{\log M_1(n)}{\lambda(n)} \; \le \; \frac14.$$
It is clear that $M_2(n) \le M_1(n)$, so to finish the proof it is enough to obtain a lower bound for $M_2(n).$

 Now, for $x \to \infty$ and $\frac{x^\rho}{\log x} \le q \le x^\rho$ (with $0 < \rho < \frac12$) choose
$$\Cal P = \Cal P(x, q) = \big\{p \le x \; \big | \; p \equiv 1 \hskip -4pt\pmod{q}\big\},$$
to be a Bombieri set relative to $x$ where $q$ is a prime number
(Bombieri prime). By Lemma
2.4 we have the
asymptotic relation
$\#\Cal P(x, q)
\sim
\frac{x}{\phi(q) \log x}.$  In order to satisfy the condition
$\prod\limits_{p\in \Cal P} p \; \le n$, we choose $x \sim q \log n.$ Without loss of generality, we may choose
$q = x^{\rho}$ for some $0 < \rho < \frac12.$ It follows that $$x^{1-\rho} \sim \log n, \qquad \log x \sim
\frac{\log\log n}{1-\rho}, \qquad \#\Cal P = \#\Cal P(x, q) \sim \frac{x}{\phi(q) \log x} \sim \frac{(1-\rho)\log
n}{\log\log n}.$$
Consequently
$$\prod_{p, p'\in \Cal P} \text{gcd}(p-1,\; p'-1) \; \ge \; q^{(\#\Cal P)^2}
\; \ge \; \left(x^\rho\right)^{\frac{(1-\rho)^2(\log n)^2}{(\log\log n)^2}   }
\sim e^{\frac{\rho(1-\rho)(\log n)^2}{\log\log n}}.$$
Let now $\rho$ go to $\frac 12$ and the theorem is proved. \qed

\enddemo


\vskip .20in

\Refs\nofrills{References}
\refstyle{A}
\widestnumber\key{MMMi}
 \vskip .10in




\ref
\key Bo
\by E. Bombieri
\paper On the large sieve
\jour Mathematika
\vol 12
\yr 1965
\pages 201-225
\endref

\vskip 1pt

\ref
\key Bu
\by L.M. Butler
\paper A unimodality result in the enumeration of
subgroups of a finite abelian group
\jour Proc. Amer. Math. Soc.
\vol 101
\yr 1987
\pages 771-775
\endref

\vskip 1pt

\ref
\key CF
\by J.W.S. Cassels, A. Fr\"ohlich
\book Algebraic Number Theory, Thompson Book Company
\yr 1967
\pages 218-230
\endref

\vskip 1pt




\ref
\key Da
\by H. Davenport
\book Multiplicative Number Theory
, Springer-Verlag, GTM {\bf 74}
\yr 1980
\endref

\vskip 1pt

\ref
\key De
\by J.B.  Dennin
\paper The genus of subfields of $K(n)$
\jour Proc. Amer. Math. Soc.
\vol 51
\yr 1975
\pages 282-288
\endref

\vskip 1pt
\ref
\key DDMS
\by        J.D. Dixon M.P.F. du Sautoy, A. Mann, D. Segal
\book Analytic Pro-p-Groups,
Cambridge University Press
London Math. Soc. Lecture Note Series {\bf 157}
\yr 1991
\endref

\vskip 1pt

\ref
\key FJ
\by M.D. Fried, M. Jarden
\book Field Arithmetic, Springer-Verlag
\yr 1986
\endref


\vskip 1pt

\ref
\key GLNP
\by  D. Goldfeld, A. Lubotzky, N. Nikolov, L. Pyber
\paper Counting primes, groups and manifolds
\jour to appear:  Proc. Nat. Acad. Sci.,  U.S.A
\vol
\yr 
\pages
\endref




\vskip 1pt

\ref
\key Kl
\by B. Klopsch
\paper Linear bounds for the degree of subgroup growth in terms of the
Hirsch length
\jour Bull. London. Math. Soc.
\vol 32
\yr 2000
\pages 403-408
\endref


\vskip 1pt

\ref
\key La
\by S. Lang
\book Introduction to Modular Forms,
Springer-Verlag
\yr 1976
\endref

\vskip 1pt

\ref
\key Lan
\by E. Landau
\book Primzahlen, Chelsea Publishing Company
\yr 1953
\endref


\vskip 1pt

\ref
\key Li1
\by U.V. Linnik
\paper On the least prime in an arithmetic progression I. The basic theorem
\jour  Rec. Math. [Mat. Sbornik] N.S.
\vol 15(57)
\yr 1944
\pages 139-178
\endref

\vskip 1pt

\ref
\key Li2
\by U.V. Linnik
\paper On the least prime in an arithmetic progression. II. The
Deuring-Heilbronn\break  phenomenon
\jour   Rec. Math. [Mat. Sbornik] N.S.
\vol 15(57)
\yr 1944
\pages 347-368
\endref





\vskip 1pt
\ref
\key LN
\by A. Lubotzky, N. Nikolov
\paper Subgroup growth of lattices in semisimple Lie groups
\jour to appear
\vol 
\yr 
\pages 
\endref

\vskip 1pt

\ref
\key LS
\by A. Lubotzky, D. Segal
\book Subgroup growth, Progress in Mathematics, Birkhauser
\yr 2003
\endref

\vskip 1pt

\ref
\key Lu
\by A. Lubotzky
\paper Subgroup growth and congruence subgroups
\jour Invent. Math.
\vol 119
\yr 1995
\pages 267-295
\endref


\vskip 1pt


\ref
\key MM
\by M. Ram Murty, V. Kumar Murty
\paper A variant of the Bombieri-Vinogradov theorem
\jour Canadian Math. Soc. Conference Proceedings
\vol  7
\yr 1987
\pages 243-272
\endref

\vskip 1pt

\ref
\key MM
\by M. Ram Murty, V. Kumar Murty, N. Saradha
\paper Modular forms and the Chebotarev density theorem
\jour Amer. J. Math.
\vol  110, no. 2
\yr 1988
\pages 253--281
\endref









\vskip 1pt

\ref
\key Pe
\by H. Petersson
\paper Konstruktionsprinzipien f\" ur Untergruppen der Modulgruppe
mit einer oder zwei Spitzenklassen
\jour J. Reine Angew. Math.
\yr 1974
\pages 94-109
\endref

\vskip 1pt

\ref
\key PR
\by V. Platonov, A. Rapinchuk
\book Algebraic Groups and Number Theory, Academic Press
\yr 1991
\endref

\vskip 1pt



\ref
\key Ra
\by S. Ramanujan
\paper Highly composite numbers
\jour Proc. London Math. Soc. (2)
\vol XIV
\yr 1915
\pages 347--409
\endref

\vskip 1pt

\ref
\key Sh
\by  A. Shalev
\paper On almost fixed
point free automorphisms
\jour J. Algebra
\vol 157
\yr 1993
\pages 271-282
\endref



\vskip 1pt

\ref
\key Su
\by M.Suzuki
\book Group Theory 1 ,Springer-Verlag, Grundlehren Math. Viss. {\bf 247}
\yr 1982
\endref




\vskip 1pt

\ref
\key Vi
\by  A.I. Vinogradov
\paper On the density conjecture for Dirichlet L-series
\jour Izv. Akad. Nauk SSSR Ser. Mat.
\vol 29
\yr 1965
\pages 903-934
\endref
\vskip 1pt

\ref
\key W
\by  A. Weil
\paper Sur les "formules explicites" de la thŽorie des nombres premiers
\jour Comm. SŽm. Math. Univ. Lund [Medd. Lunds Univ. Mat. Sem.] 
\vol 
\yr 1952
\pages 252-265
\endref





\endRefs
\enddocument